\documentclass[10pt]{article}

\input epsf
\usepackage{enumerate, fullpage}
\usepackage{array}
\usepackage{graphicx}
\usepackage{color}
\usepackage{amsfonts}
\usepackage{amsmath}
\usepackage{url}
\usepackage{multimedia}
\usepackage{ulem}
\usepackage{subfig}

\usepackage{booktabs}
\usepackage{array, supertabular}
\usepackage{multirow}
\usepackage{longtable}

\def\b1{\mathbf{1}}

\newcommand{\beq}{\begin{equation}}
\newcommand{\eeq}{\end{equation}}

\newcommand{\beqnr}{\begin{eqnarray}}
\newcommand{\eeqnr}{\end{eqnarray}}

\newcommand{\benum}{\begin{enumerate}}
\newcommand{\eenum}{\end{enumerate}}
\newcommand{\argmax}{\mathop{\rm argmax}}

\newcommand{\QED}{\rule{.1in}{.1in}}

\newcommand{\cH}{{\cal H}}

\newcommand{\cP}{{\cal P}}
\newcommand{\cQ}{{\cal Q}}

\newtheorem{DE}{Definition}[section]

\newtheorem{LE}[DE]{Lemma}

\newtheorem{CO}[DE]{Corollary}

\newcommand{\qed}{\mbox{}\hspace*{\fill}\nolinebreak\mbox{$\rule{0.7em}{0.7em}$}}

\advance \topmargin by -\headheight
\advance \topmargin by -\headsep

\evensidemargin \oddsidemargin

\begin{document}

\begin{center}

\Large {\bf Models for managing the impact of an epidemic}
\vskip 0.1cm
\vskip 0.5cm

\Large

Daniel Bienstock and A. Cecilia Zenteno \\
Columbia University \\
New York\\
(Preliminary version)\\
February, 2012\\

\vspace{.1in}
\small version Sat.Mar.10.153610.2012
\normalsize
\end{center}

\section{Introduction}\label{sec:intro}
In this paper we consider robust models for emergency staff deployment in the event of a flu pandemic.  We focus on managing critical
staff levels at organizations that must remain operational during such an event,  and develop methodologies for managing emergency resources with the goal of minimizing the impact of the pandemic. We present numerical experiments using realistic data to
study the effectiveness of our approach.

A serious flu epidemic or pandemic, particularly one characterized by
high contagion rate, would have extremely damaging impact on a large,
dense population center.  The 1918 influenza pandemic is often seen as a worst-case scenario as it arguably represents the most devastating pandemic in recent history, having killed more than
20 million people worldwide \cite{chowell2007comparative, mills2004transmissibility, taubenberger20061918}.   However, even a much milder epidemic would
have vast social impact as services such as health care, police and utilities became severely hampered by staff shortages. Workplace absenteeism might also become a serious concern \cite{Ch6_protectstaffpandemic}; for example, the New Zealand government has predicted overall absenteeism levels as high as 40\% (\cite{lee2007effectiveness} and references therein)\footnote{In the healthcare setting, the opposite may take place: health care workers reportedly avoid calling in sick during an emergency \cite{convYoko, bridges2003transmission}}

In this study we focus on managing the
inevitable staff shortfall that will take place in the case of a severe
epidemic.  We take the viewpoint of an organization that
seeks to diminish decreased performance in its operations as the epidemic
unfolds, by appropriately deploying available resources, but which is not directly attempting to control the number of
people that become infected. In contrast to our focus, much
valuable research has been directed at addressing the
epidemic itself;  such work studies public health measures that would reduce the epidemic's severity and its direct impact,
for example by managing the supply of vaccines and antivirals
(see \cite{longini2004contWantiviral, hill2003critical, longini2005contAtsource}).  While we do not address this topic, it seems plausible that
our methodology for handling robustness will apply in this setting, as well.

We are interested in infrastructure of critical social value, such
as hospitals, police departments, power plants and supply chains, among many examples of
entities that must remain operative even as staff levels become low. In
cases such as police departments, staff would
likely be more exposed to the
epidemic than the general public and (particularly if vaccines are in short
supply, or apply to the wrong virus mutation) shortfalls may take place just when there is greatest demand for services.
Power plants and waste water treatment
plants are examples of facilities whose operation
will be degraded as staff falls short and which
probably require minimum staff levels to operate at all \cite{hoffbuhr2006utilities}.  Supply chains
would very likely be significantly slowed down as their staff is depleted,
resulting, for example, in food shortages \cite{osterholm2005preparing}.
In all these cases, organizations cannot implement ``work from home'' strategies as urged for private businesses by the Centers of Disease Control and Prevention (CDC) and the United States Department of Labor \footnote{\url{http://www.osha.gov/Publications/employers-protect-workers-flu-factsheet.html}}.


A pandemic contingency plan for a large organization (such as a city government)
would include resorting to emergency (or ``surge'') sources for additional staff: for example
by temporarily relying on personnel from outlying, less dense, communities.  Such additional
resources are likely to be significantly constrained in quantity, duration and
rate, among other factors.  Such emergency staff deployment plans would
entail some complexity in design, calibration and implementation, but as a result of other disasters such, as Hurricane Katrina in 2005 and the anthrax attacks in 2001, it is now agreed that there is a compelling need for emergency planning \cite{gershon2007homeHCchallenges}.


From a planner's perspective, the task of managing future resource levels
during an epidemic is complex, partly because of uncertainties regarding the behavior of the epidemic, in particular, uncertainty in the contagion rate.  The evolution of the
contagion rate is a function of poorly understood
dynamics in the mutations of the different strains of the flu virus and environmental
agents such as weather \cite{chowell2007comparative, flu_coldweather}.  In addition to uncertainty, a decision-maker will also likely be constrained by logistics. In particular, it may prove impossible to carry out large changes to staff deployment plans
on short notice, particularly if such staff is also in demand by other
organizations (as might be the case during a severe epidemic).  We will return
to these issues in Sections \ref{subsec:deploymentstrat} and \ref{thefuture}. As a consequence of the two factors (uncertainty, and logistical constraints) a decision-maker may commit
too few or too many resources - in this case emergency staff- or perhaps at the wrong time,
if there turns out to be
a mismatch between the anticipated level of staff shortage and what actually transpires (after the resources have been committed).

We present models and methodology
for developing emergency staff deployment levels which optimally hedge
against the uncertainty in the evolution an the epidemic while accounting for
operational constraints.
Our approach overlays adversarial models to describe the contagion rate on the classical SEIR model for describing
epidemics. The resulting robust optimization models are nonconvex and
large-scale; we present convex approximations and
algorithms that prove numerically accurate and efficient, and we
study their behavior and the policies they produce under a range of scenarios.

This paper is organized as follows.  Section \ref{sec:fluprep} addresses issues related to surge capacity planning; Section \ref{sec:SEIR} describes the classical (non-robust) SEIR model; Sections \ref{sec:umeasures},
\ref{sec:SEIRrobust} and \ref{sec:robust} contain the description of our initial robust model,  Section \ref{sec:experiments} describes our experiments, and our algorithms are detailed in Appendix \ref{appendix}.

\section{Planning for workforce shortfall}\label{sec:fluprep}

An influenza pandemic would severely stress
the operational continuity of social and business structures
through staff shortages.   Staff shortfall directly resulting from individuals becoming sick could be intensified by policy or absenteeism. For example, during the last H1N1 influenza outbreak, CDC recommended that people with influenza-like symptoms remain at home until at least 24 hours after they are or appear to be free of fever. In the particular case of health care workers it is advised that they refrain from work for at least 7 days after symptoms first appear  (see \cite{url:cdc} for additional details). Moreover, staff shortages may occur not only due to actual illness, but also
from illness among family members, quarantines, school closures (combined with lack of child care),
public transportation disruptions, low morale, or because workers could be summoned to comply with public service obligations \cite{Ch6_protectstaffpandemic, PandPrep_Gomersall07}. Indeed, employees who have been exposed to the disease
(especially those coming into contact with an ill person at home) may also be asked to stay at home and monitor their own symptoms.

 Altogether, the direct and indirect staff shortfall caused by an epidemic,
in a worst-case scenario, could give rise to an extended time period during which 20 to 40\% of the workforce will be absent, public health and utility professionals predict \cite{hoffbuhr2006utilities}.
Even though the outlook would be dire, organizations that provide critical infrastructure services such as health care, utilities,
transportation,
and telecommunications,  should clearly continue operations and are exhorted to plan accordingly \cite{url:flu.gov}. (See
\cite{morse2003building} for additional background on emergency staff
planning.)

The Department of Health and Human Services (HHS) and the CDC have provided guidelines to help businesses and their employees in planning for an
influenza outbreak. In addition to recommendations addressing
the spread of disease and antiviral drug stockpiling, there is a focus on staff planning, which is
the subject of our study \cite{url:flu.gov, url:cdc}.
Additionally, there are federal and state programs such as the New York Medical Reserve Corps whose mission is to organize volunteer networks prepare for and respond to public health emergencies, among other duties \cite{url:NYCMedResCorps}.

HHS and CDC urge organizations to identify critical staff requirements needed to maintain operations during a pandemic. In particular, organizations should develop detailed criteria to determine when to trigger the implementation of an emergency staffing plan. Most significantly, organizations should identify the minimum number of staff
needed to perform vital operations.  For example, in the case of water treatment plants, approximately 90\% of the personnel is critical for keeping the utility running; for refineries, losing 30\% of their staff would force a shutdown \cite{hoffbuhr2006utilities}.

In the particular case of hospitals, an effective \textit{contingency staffing}
plan should incorporate information from health departments and emergency management authorities at all levels, and  would build a data base for alternative staffing sources
(e.g., medical students). For additional details, see \cite{url:flu.gov, Ch2_SurgeCapacity}. In New York City, for example, after the events of 9/11, Columbia University  created a database of volunteers to be recruited and trained both in basic emergency preparedness and their disaster functional roles \cite{morse2003building}. At the city level the NYC Medical Reserve Corps ensures that a group of health professionals ranging from physicians to social workers is ready to respond to health emergencies. The group is pre-identified, pre-credentialed, and pre-trained to be better prepared in the wake of a crisis. Similar emergency staff backup plans
could be implemented in all other cases of utilities and social infrastructure \cite{url:NYCMedResCorps}.

In spite of all these efforts, it seems clear that much remains to be done and that a severe epidemic would
place extreme strain on infrastructure. A good example is provided by the 2009 Swine Flu epidemic.  Even
though the virus mutation caused few fatalities, and a successful vaccine became available, New York hospitals
were severely stressed \cite{url:NYT_H1N1flu}: ``The outbreak highlighted many national weaknesses: old, slow vaccine technology; too much reliance on foreign vaccine factories; {\bf some major hospitals pushed to their limits by a relatively mild epidemic}'' (our emphasis).

From our perspective, the uncertainty concerning the timeline and severity of the pandemic brings substantial complexity to the problem of deploying replacement staff.  This problem,
which will be the core issue that we address, is relevant because significant preplanning must take
place and it is unlikely that major quantities of
additional workforce can be summoned on a day-by-day basis.

\subsection{Declaring an epidemic}\label{subsubsec:declare}  A technical point that we will return to below concerns when, precisely, an epidemic is ``declared''. In the case of an infectious disease such as influenza,
an initially slow accumulation of cases followed by a more rapid increase in incidence \cite{whenis1epidemic1epidemic} is viewed by
epidemiologists as an epidemic. However, this definition is too general for
planning purposes. This is a significant issue since emergency action plans would be activated when the epidemic is declared.

In the United States, the CDC declares an influenza epidemic when death rates from pneumonia and influenza exceed a certain threshold \cite{url:CDCFluWeekly}. Each week, vital statistics of 122 cities report the total number of death certificates received and the proportion of which are listed to be due to pneumonia or influenza. This percentage is compared against a seasonal baseline, which in turn was computed using a regression model based on historic data. There is a different baseline for each week of the year to capture the different seasonal patterns of influenza-like illnesses (ILI) (the epidemic threshold sits 1.645 standard deviations above the seasonal baseline). This type of measure is not specific for the United States. In \cite{lee2007effectiveness}, the authors' base case assumed that the duration of an influenza pandemic in Singapore was defined as the period during which incident symptomatic cases exceeded 10\% of baseline ILI cases.

Motivated by this discussion, we will use the convention that an epidemic is known to be present as soon as the number of (new) infected individuals on a given time period exceeds a small percentage of the overall population, e.g. $0.93\%$, which corresponds to the national epidemic threshold of $6.5\%$ for week $40$\cite{url:CDCFluWeekly}.

From the point of view of an organization, of course, action need not
wait until an ``official'' epidemic declaration and would instead rely on its
own guidelines to possibly implement a preparedness plan at an earlier point.
However, we expect that the
mechanism underlying declaration will be the similar.


\section{Modeling influenza}\label{sec:SEIR}

Influenza is an acute, highly contagious respiratory disease caused by a number of different virus strains. While most people recover within one or two weeks without requiring medical treatment, influenza may cause lethal complications such as pneumonia. For certain virus strains this is especially pronounced in the case of susceptible groups such as young children, the elderly or people with certain medical preconditions.

Yearly influenza outbreaks occur because the different virus strains acquire adaptive immunity to old vaccines by frequent mutations of their genetic content \cite{earn2002ecologyflu}. According to the CDC, seasonal influenza causes an average of 23,600 deaths and more than 200,000 hospitalizations in the U.S. \cite{url:flu.gov}. When a larger genomic mutation occurs the virus derives into a new one for which humans have little or no immunity resulting in a pandemic: the disease may rapidly spread worldwide, possibly with high mortality rates. 

There were three flu pandemics in the twentieth century, the worst of which
occurred in 1918; known as the ``Spanish flu'', it killed 20-40 million people worldwide. Milder pandemics
occurred in 1957 and 1968. The most recent pandemic occurred in 2009 - 2010 with the surge of the H1N1 virus. Current fears that the avian strain could develop into another pandemic have caused interest in modeling the spread of influenza and evaluating possible emergency management strategies.

The critical difficulty in modeling the impact of a future influenza
pandemic
is our inability to accurately predict the spread of disease on a given population.   In this section we provide a description of the (nonrobust) model which we will modify in order to follow the evolution of the disease and its spread among members of a given population. We focus on
the progress of the epidemic within a workforce group of interest. We then discuss weaknesses of the model and the use of robust methods to remedy these shortcomings.


\subsection{SEIR -- a basic influenza model}\label{subsec:compartmodel}

A classic simple model for influenza is the so called $SEIR$ compartmental model. These kind of models describe, in a deterministic fashion, the evolution of the disease by partitioning the population into groups according to the progression of the disease. The terminology $SEIR$ describes the transition of individuals between compartments: from susceptible ($S$) to being in an incubation period ($E$) to becoming infectious ($I$) and finally to the removed class ($R$). The latter could include death due to infection as well as recovered people, depending on the model.

We make the following standard assumptions \cite{brauer2008mathepi}:

\begin{enumerate}
    \item There is a small number $I_0$ of initial infectives relative to the size of the total population, which we denote by $N$.
    \item The rate at which individuals become infected is given by the product of the probability that at time $t$ a contact is made with an infectious person, $\beta_t$, the average constant social contact rate $\lambda$, and the likelihood of infection, $p$, given that a social contact with an infectious person
has taken place. The probability $\beta_t$ changes with time because we assume it depends on the number of infectious agents in the population.
    \item The latent period coincides with the incubation period; that is, exposed individuals are not infectious. We assume they proceed to the $I$ compartment with rate $\mu_{E}$.
    \item Infectives leave the infectious compartment at rate $\mu_{RR}$.
    \item We assume there is only one epidemic wave. Thus, people who recover are conferred immunity.
    \item The fraction of members that do not die from the disease, when removed from the infectious class, is given by $0 \leq f<1$.
    \item We do not include births and natural deaths because influenza epidemics usually last few months. We also omit any migration. In other words, excluding deaths by disease, the total population remains constant.
\end{enumerate}

SEIR models are usually described as a system of nonlinear ordinary differential equations (see for example \cite{brauer2008mathepi, NewallCEAaustralia}). For our purposes, we use a discrete-time Markov chain type approximation along the lines of \cite{larson2007infprogressmodel} and similar to those found in \cite{allen2000comparison, allen1991discrete, allen2008basic}:

\begin{eqnarray}
	   S_{t+1} &=& S_t(e^{-\lambda \beta_t p}) \nonumber \\
       E_{t+1} &=& E_t(e^{-\mu_E}) + S_t(1 - e^{-\lambda \beta_t p})\\       \label{eq:SEIRsystem}
	   I_{t+1} &=& I_t f (e^{-\mu_{RR}}) + E_t(1-e^{-\mu_E}) \nonumber \\
	   R_{t+1} &=& R_t + I_t (1-e^{-\mu_{RR}}). \nonumber
\end{eqnarray}


Compartmental models incorporate a number of assumptions to describe social contact dynamics. First, the number of social contacts with infectious people for an arbitrary person is thought of as a Poisson random variable with rate $\lambda \beta_t$; we use one day as time unit. Second, the models assume homogeneous mixing, that is, all individuals have a fixed average number of contact rates per unit of time and are all equally likely to meet each other. Thus, the probability that a contact is made with an infectious person, $\beta_t$ is given by $I_t / N$. The daily infectious contact rate $\lambda p$ is usually written as $\lambda$ \cite{brauer2008mathepi, allen1991discrete} and taken as a constant throughout the epidemic. We will remove this assumption later when we make the transmissibility parameter $p$ explicit.

If the number of initial infectives is relatively very small compared to the whole population ($S_0 \sim N$), then a newly introduced contagious individual is expected to infect people at the rate $\lambda p$ during the expected infectious period $1/\mu_{RR}$. Thus, each initial infective individual is expected to transmit the disease to an average of

\begin{equation}\label{eq:R_0_1}
    R_0 = \frac{\lambda \, p} {\mu_{RR}}
\end{equation}
individuals. $R_0$ is called the \textit{basic reproduction number} (also known as \textit{ basic reproduction ratio} or \textit{ basic reproductive rate}.) It is without doubt the most important quantity epidemiologists consider when analyzing the behavior of infectious diseases \cite{brauer2008mathepi}. Its relevance derives mainly from its threshold property: when $R_0 < 1,$ the disease does not spread fast enough and there is no epidemic; when an epidemic does take place, it is of interest to see how much bigger than 1 $R_0$ is. Because we assume that an epidemic will take place \textit{a priori}, $R_0$ is not of particular interest to us; however, it is presented here for completeness.

\subsubsection{Keeping track of staff availability}\label{subsubsec:SEIR_2ndSet}

We are interested in tracking workforce
availability at a particular organization during the epidemic; following previous work \cite{arino2006simple, gardam2007prophylhcw} we divide the population into two groups: (1) the general population and (2) the workforce under consideration. Individuals from the latter group could have a very different exposure to the epidemic. For example, people working at a water plant could have lower contact rates than average
(by virtue of having contact with few individuals during the workday) while the staff at a health clinic may not only have higher contact rates, but may also have easier access to antiviral medicines that reduce their infectiousness and the length of the infectious period.

For ease of notation, we use superscript $1$ to refer to the general population and $2$ to refer to the group of workers of interest. For $j = 1,2$, define
compartments $S^j, E^j, I^j, R^j$ corresponding to group $j$. Following the
discussion above, we allow the groups to have different contact, incubation, and recovery rates. The probability that a random contact is one with an infected person at time $t$, $\beta_t$, is now defined as

\begin{equation}\label{eq:beta}
\beta_t = \frac{\lambda^TI_t}{\lambda^TN_t},
\end{equation}

\noindent where $I_t = [I^1_t,I^2_t]$ is the vector of infectious individuals
at time $t$, $\lambda^T = [\lambda^1, \lambda^2]$ is the vector of contact rates, and $N_t = [N_t^1, N_t^2]$ denotes the size of each group at time $t$. We note that $\beta_t$ is constant across groups. We now have two parallel thinned Poisson process approximations, each with rate $(\lambda_j\beta_t p)$. The set of equations that correspond to group $j$ ($= 1,2$) is

\begin{eqnarray}
	   S^j_{t+1} &=& S^j_t e^{-\lambda_j \beta_t p} \nonumber \\
       E^j_{t+1} &=& E^j_t e^{-\mu_E} + S^j_t(1 - e^{-\lambda_j \beta_t p})\\ \label{eq:SEIRsystem_2}
	   I^j_{t+1} &=& I^j_t f e^{-\mu_{RR}} + E^j_t(1 - e^{-\mu_E}) \nonumber  \\
	   R^j_{t+1} &=& R^j_t + I^j_t (1 - e^{-\mu_{RR}}). \nonumber
\end{eqnarray}

\noindent We now give an expression for $R_0$ \cite{brauer2008mathepi, dieckmann2000mathematical}. First we note that the average contact rate for subgroup $j$ at time $0$ is given by
$$ \frac{\lambda_j N^j}{\lambda_1 N^1 + \lambda_2 N^2}.$$
Thus, the average number of new infections caused by a newly infective person introduced into an otherwise susceptible population is given by
\begin{align}
\nonumber    R_0 &= \left[ \lambda_1 \frac{\lambda_1 N^1}{\lambda_1 N^1 + \lambda_2 N^2} + \lambda_2 \frac{\lambda_2 N^2}{\lambda_1 N^1 + \lambda_2 N^2} \right]\frac{p}{\mu_{RR}}\\
        &= \frac{\lambda_1^2 N^1 + \lambda_2^2 N^2}{\lambda_1 N^1 + \lambda_2 N^2} \cdot \frac{p}{\mu_{RR}}. \label{eq:R_0_2}
\end{align}

We note that when $\lambda_1 = \lambda_2$ (\ref{eq:R_0_2}) reduces to (\ref{eq:R_0_1}).

\subsubsection{Nonhomogeneous mixing and social distancing}\label{subsubsec:NonHomMix_SocDist}

We initially assumed that the population mixed homogeneously, that is the contact rates $\lambda_1$ and $\lambda_2$ remained constant throughout. However, the homogeneous mixing and mass action incidence assumptions have clear pitfalls; individuals from each compartment are hardly indistinguishable in terms of their social patterns and likeliness of infection. As more people acquires the disease, members in the population may become more careful of their social contacts and infectious people may stay at home. Thus, following \cite{allen1991discrete} we make the assumption that

\begin{equation}\label{eq:nonhomlambda}
    \lambda^j_t = \Lambda_j \frac{S^j_t + E^j_t + R^j_t}{N^j_t}, \quad j = 1,2,
\end{equation}
where $\Lambda_j$ are fixed constants, $j = 1,2.$ Using this definition, $\lambda^j_t$ decreases whenever there is a high number of infectious agents in the population. As mentioned in \cite{allen1991discrete}, other functional forms are possible; we rely on (\ref{eq:nonhomlambda}) because it provides a
simple way to capture changes in contact rates as a function of  severity of the epidemic.

Additionally, we also consider a scenario in which authorities impose social distancing measures as soon as the epidemic is declared. A similar situation took place in
Mexico during the last 2009 H1N1 epidemic, where venues such as schools, movie theaters, and restaurants were forced to close temporarily. It is estimated that the transmission of the disease was diminished by $29$\% to $37$\% \cite{chowell2011H1N1Mexico}. We incorporate this element by multiplying the contact rates by an additional dampening factor when the epidemic is considered declared and until the rate of growth of daily infectives is below some threshold (we refer the reader to section \ref{subsubsec:declare}). Effectively, the contact rate for group $j$ at time $t$ becomes $\lambda^j_t = \theta \Lambda_j \left( \frac{S^j_t + E^j_t + R^j_t}{N^j_t} \right)$ ($0 < \theta < 1$) when the epidemic is officially ongoing; otherwise, it remains as per equation (\ref{eq:nonhomlambda}).

\section{Robustness in planning}\label{sec:SEIRrobust}

In planning the response to a future or impending epidemic, one would need
to rely at least partly on an epidemiological model, and such a model
would have to undergo careful calibration in order to be put
to practical use. This is especially the case for the SEIR model presented above, which is rich in parameters that need to be estimated
to fully define the trajectory of the epidemic.
In the case of a flu pandemic caused by an unknown virus strain,
new parameter values would need to be promptly estimated as the epidemic emerges.  Ideally,
robust statistical inference would provide information on all parameters,
though data paucity would present a challenge.

On the positive side, the infectious and incubation periods
can usually be independently estimated via clinical monitoring of infected
agents, either by observation of transmission events or by the use of more
detailed techniques \cite{keeling2008modeling}.

On the other hand, it is not clear how to accurately estimate the\textit{ transmission rate} $\lambda_j p$. One approach is to approximate the basic reproductive ratio $R_0$ and the mean infectious period, and then use equation (\ref{eq:R_0_1}). There are multiple ways of estimating $R_0$; see, for example, \cite{chowell2007comparative}. Given that the definition of $R_0$ is based on the early stages of the epidemic, one should examine its early behavior. However, the progression of the epidemic in its initial stages could fluctuate widely because of the very small number of initial cases, making the fitting process more difficult \cite{keeling2008modeling}.

Direct estimation of the transmission rate gives rise to a
number of challenges \cite{wallinga2006using, wu2006reducing}. First, pandemic influenza
(along with smallpox and pneumonic plague)
has not been present in modern times frequently enough  so
as to gather sufficient data for accurate estimation.
Second, for existing age-specific
transmission models, there are more parameters than observations on risk
of infection for each age class.
The infectivity of the disease can be approximated using serologic data and contact rates are usually estimated from census and transportation data after making assumptions about the contact processes that reduce the number of unknowns to the number of age classes.  However, both contact or infection rates estimates
can serve only as a baseline; a population could change its behavior significantly during a severe epidemic (due to school closures, for example); further,
environmental changes (e.g. weather changes) could also have a significant impact on the virus transmissibility \cite{flu_coldweather}.

The infectivity parameter
$p$ is particularly problematic because (at least partly) it reflects
the characteristics of the virus; it is usually estimated after the epidemic has taken place.   It seems very difficult to predict the evolution of new virus strains.  In fact, as far
as we know \cite{convRaulRabadan, ConvStephenMorse}) research that relates
mutations of influenza virus to infectivity is inconclusive. Previous work has proposed different upper and lower bound values for $p,$ depending, among other factors, on the geographical location of the study and the pandemic wave of interest. For example, \cite{wallinga2006using} uses the interval $[0.025, 0.5]$  while $[0.02, 0.16]$ is used in \cite{wallinga2010optimizing}. Both studies use these intervals to conduct sensitivity analysis. In another study, Larson \cite{larson2007infprogressmodel} classifies the population into three groups according to social activity levels; the three groups have $p$ value $0.07, \, 0.09$ and $0.12$ and corresponding $\lambda$ values $50, \,10$ and $2$, respectively. In summary, the precise estimation of $p$ values appears quite challenging, especially prior to or even during an epidemic.


\begin{figure}[h!]
\centering
\subfloat[An absolute perspective in time: Curves start from the moment infectives are introduced into the susceptible population.]{
\includegraphics[width=0.45\textwidth]{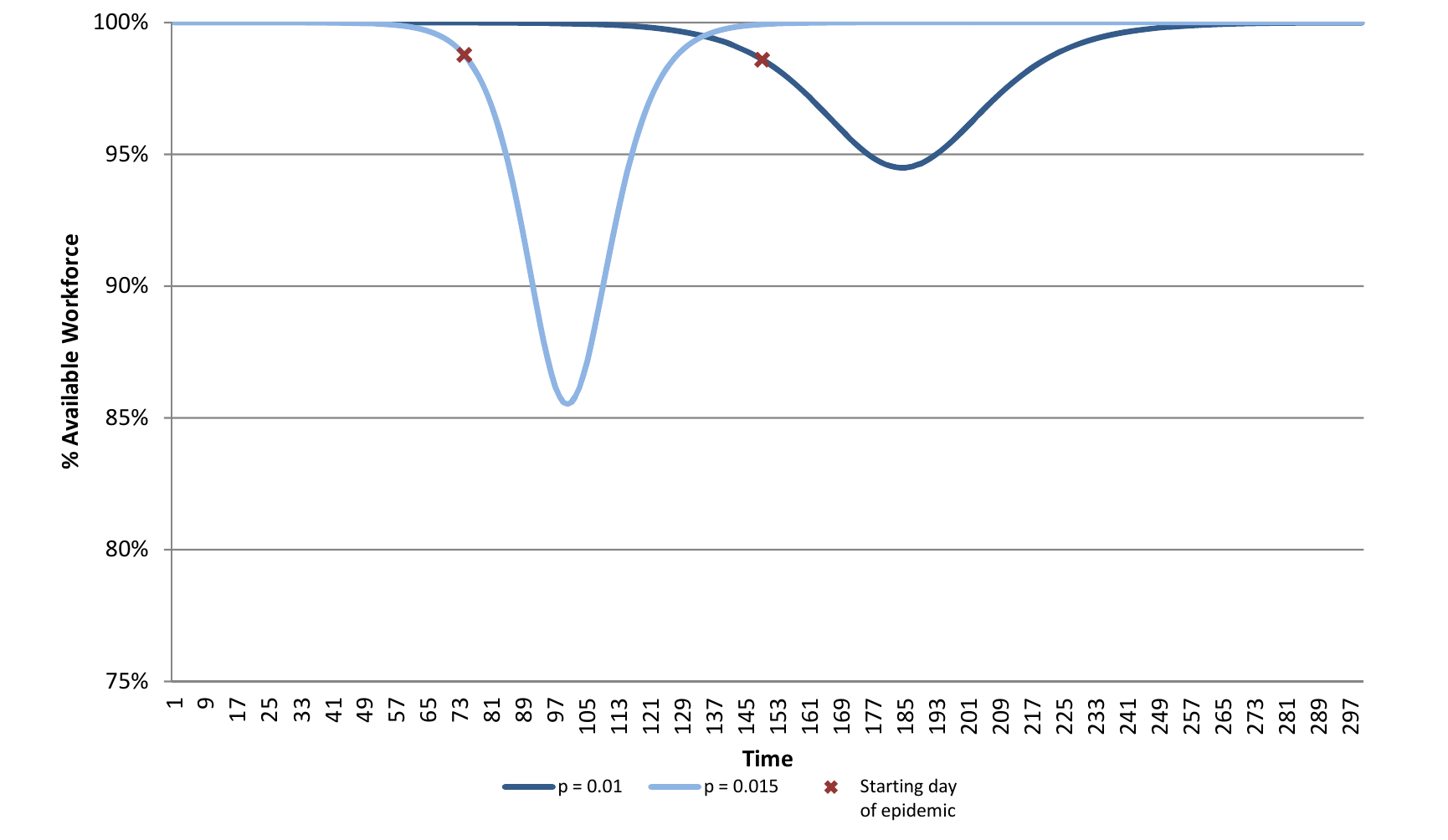}
\label{chart:subAvWorkers01a}
}
\hspace{10pt}
\subfloat[A planner's perspective: Curves for \textit{both epidemics} start from the moment an epidemic is declared.]{
\includegraphics[width=0.45\textwidth]{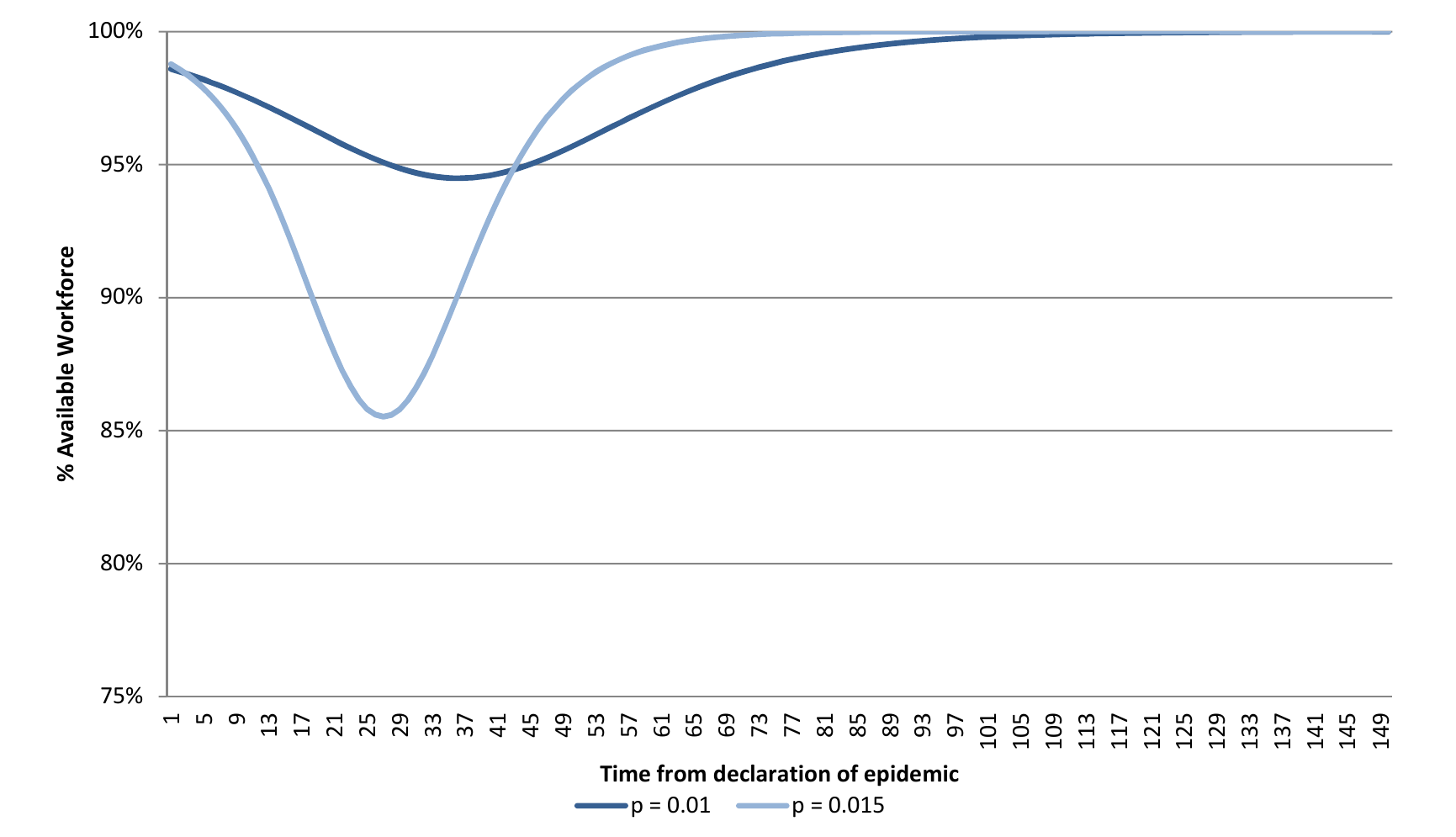}
\label{chart:subAvWorkers01b}
}
\caption{Availability of workforce as epidemic progresses for different values of $p$.}
\label{chart:AvWorkers}
\end{figure}


Given this uncertainty, it is likely that basing the response to an epidemic on a fixed estimate for
$p$ is incorrect. To illustrate the impact of such a decision we refer the reader to Figure \ref{chart:AvWorkers}. It shows the availability of staff as a function of time, for different two values of $p$, from two different perspectives.
Figure \ref{chart:subAvWorkers01a} shows how the workforce becomes ill at the actual time it happens.
In contrast, \ref{chart:subAvWorkers01b} presents these curves all starting from the time the epidemic is declared.
Indeed, a planner deploying a contingency plan at the moment an epidemic is declared would be interested in the preparing for different scenarios according to what is shown in \ref{chart:subAvWorkers01b}, rather than \ref{chart:subAvWorkers01a}. This point will be touched upon again in Section \ref{subsec:deploymentstrat}.

Consider a baseline of $90\%$, i.e. the system is considered to be performing poorly
if fewer than $90\%$ of the staff is available.  For $p = 0.015$ this period
spans days $18$ through $37$ (from the declaration of the epidemic), whereas for $p = 0.01$ the baseline is not reached. As expected, the epidemic is more severe
for higher values of $p$; however, somewhat lower values of $p$ result in
longer-lasting epidemics.  In particular, $p = 0.01$ results in an extended
period of time (days $13$ through $69$) where even though staff availability
is above the baseline, it is still significantly below $100 \%$.  If, in the above example, $p$ were unknown, a planner would have to carefully ration scarce resources over a nearly month-long period.  If the planner
were to assume a {\em fixed} value for $p$ throughout the epidemic, then the
example suggests allocating comparatively higher levels of surge staff to
earlier periods of time, to overcome the higher shortfalls to be expected in
the case of higher values for $p$.  Of course, this higher level of initial
allocation needs to be carefully chosen to obtain maximum reward.



\begin{figure}[h!]
\centering
\subfloat[An absolute perspective in time: Curves start from the moment infectives are introduced into the susceptible population.]{
\includegraphics[width=0.45\textwidth]{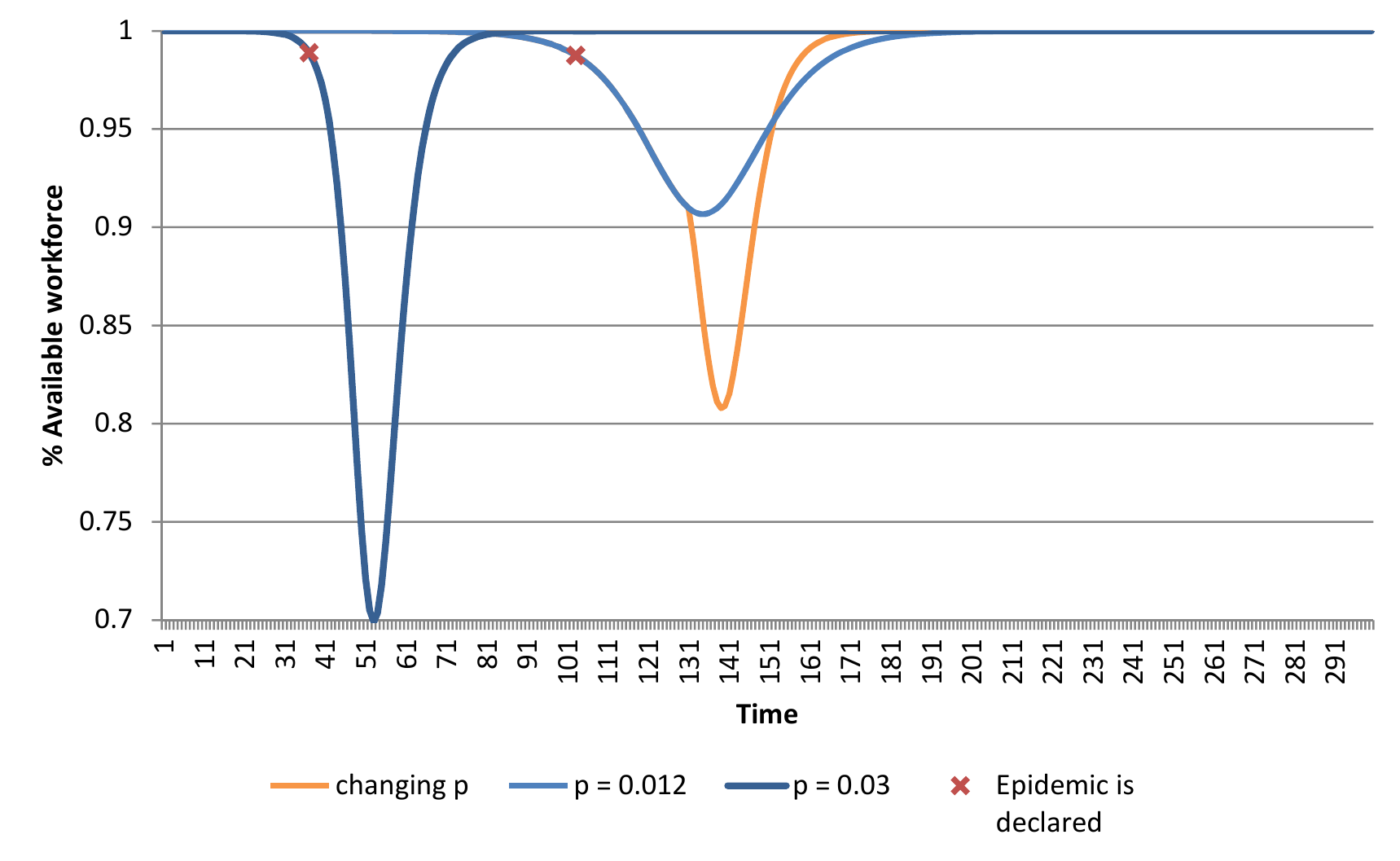}
}
\hspace{10pt}
\subfloat[A planner's perspective: Curves for \textit{all three epidemics} start from the moment an epidemic is declared.]{
\includegraphics[width=0.45\textwidth]{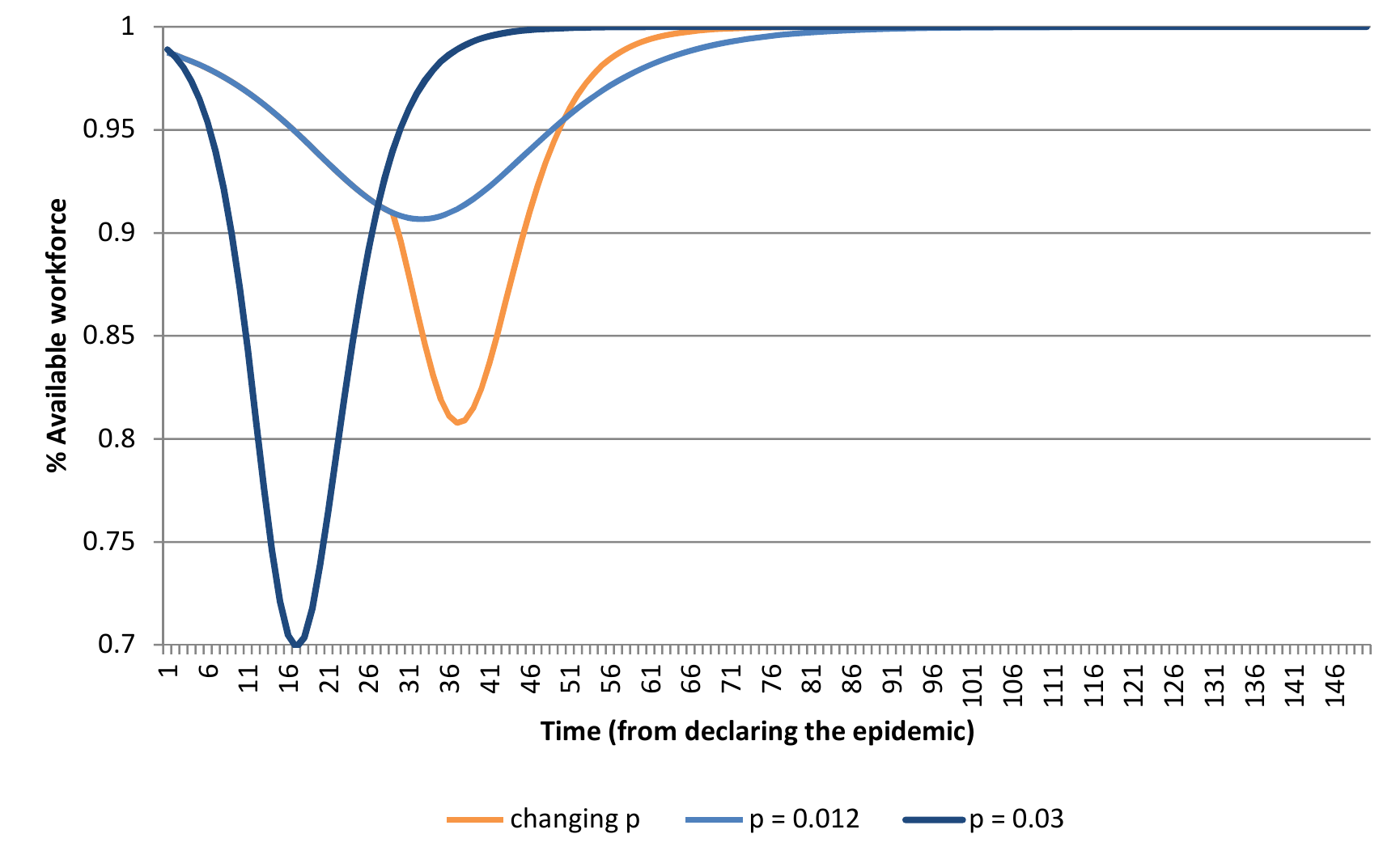}
}
\caption[Availability of workforce as epidemic progresses for different values of $p$.]{Availability of workforce as epidemic progresses for different values of $p$.}
\label{chart:pchange_03}
\end{figure}

Consider now Figure \ref{chart:pchange_03}, which shows the impact of a {\em change} in $p$ in the
midst of an epidemic. Here, $p$ changes from $0.012$ to $0.03$
on day $130,$ $28$ days after the epidemic has been declared.  Thus, for more than half of the epidemic, the disease spreads slowly and
even though the $90\%$ baseline is approached, it is not reached. However,
after the change in $p$ the epidemic becomes severe and a significant shortfall arises. This type of variability would be especially problematic when resources are limited.  The epidemic, on the basis of observations of its initial progression,
would likely be classified as relatively mild, and, perhaps,
action might be taken to at least partially
abandon the surge staff buildup. However, if the a change in $p$ as shown
in Figure \ref{chart:pchange_03} were plausible, then
a careful planner would have to hedge by hold backing
staff so as to handle the potentially critical situation in later periods.
Should the change in $p$ not take place,
the chart indicates that any held back staff would essentially be wasted.  And
should the change take place somewhat earlier, then the opportunity cost
is significantly higher.  The critical question, of course, is whether this
type of virus behavior is possible.
As far as we can tell, current knowledge of the
influenza virus cannot categorically reject a change in $p$ as shown,
although planning for such a change might be construed as an overly conservative
action. Thus, it remains to be seen if a robust strategy that can protect
against a change in $p$ entails higher cost or other compromises.


In this paper we model the variability of the product $\lambda_j p$ using
techniques derived from robust optimization.  For ease of exposition,
we assume that the infectivity probability $p$ is uncertain, and that the values for the contact rates $\lambda_j$ are fixed (and known); thus the contact probabilities $\beta_t$ are also fixed. Our goal is to produce surge staff deployment strategies that are robust with
respect to variability of $p$ in a number of models of uncertainty.   Robust optimization
provides an agnostic methodology for assessing competing allocation plans
and for computing good plans according to various criteria.  Above,
we described two such criteria.  Our optimization procedures
are fast and flexible, making it possible to evaluate multiple plans and
different levels of conservatism in a practicable time frame.

\section{Performance Measures}\label{sec:umeasures}

Personnel shortage during an epidemic may jeopardize the continuity of operations of critical infrastructure organizations. To gauge the impact of this shortfall, we consider cost functions that model two scenarios. In the first one the organization requires at least certain percentage of personnel present to operate. This could be the case of water and energy plants, as mentioned in \cite{hoffbuhr2006utilities}. The second scenario uses classic queueing theory to
measure the simultaneous impact of the decrease in number of available workers and the change in demand for the service that the organization provides. Health care institutions would be clear examples in this context.\\

\subsection{Threshold functions} \label{subsec:threshfun}

Here we model an organization that requires a minimum number of staff in order
to operate under normal conditions.  We further assume that this ``minimum'' is
a \textit{soft} constraint in the sense that if the available staff should fall
below the threshold, the organization will still manage to operate, but at
very large cost. This could be the case, for example, if the organization is
able to purchase output or hire staff from another source (such as a competitor) or if it is
able to reduce its level of service or output, at cost.

To model this behavior, we assume that operating costs at each point in time are given by a convex piecewise linear function.   This function is represented
as the maximum of $L$ linear functions, with slopes
$\sigma_L < \sigma_{L-1} < \ldots < \ldots \sigma_1 = 0$ and intercepts $k_L > K_{L-1} > \ldots > K_1 = 0$. Thus, denoting by $\omega_t$ the work force level at time $t$ and $z_t$ the cost associated with time period $t$, we have
 \begin{equation}\label{eq:costfn_threshold}
z_t = \max_{1 \le i \le L} \left\{ \sigma_i \, \omega_t + k_i \right\}.
\end{equation}

\subsection{Queueing models}\label{subsec:queue}
We are interested in modeling basic queueing systems where both the service rate
and the incoming demand rate are affected by the evolution of the epidemic.
For each time period $t$, denote the average incoming demand rate by $\zeta_t$ and let $s_t$ denote the number of available servers (staff). The system utilization at time $t$, using an M/M/$s_t$ queueing model, is given by
\begin{equation}\label{eq:rho_t}
    \rho_t = \frac{\zeta_t}{s_t \mu}.
\end{equation}
As it is well known, $\rho_t$ and related quantities are good indicators of system
performance.  In our simulations of severe epidemics,
$\rho_t$ can grow larger than $1$, a situation that depicts a system that is catastrophically saturated
(a situation observed during 2010's H1N1 pandemic\cite{url:NYT_H1N1flu, convYoko}) and consequently system performance will
drastically suffer.  Accordingly, we choose as a reasonable representative
of ``cost'' incurred at time $t$ an exponentially increasing function of $\rho_t$.
In particular we consider $e^{\rho_t - \delta}$ (for small $\delta > 0$) which we approximate
with a piecewise-linear function, much as in Section \ref{subsec:threshfun}.   Other alternatives are possible. For example, one could use a traditional measure of system performance, such
as average queue length for an M/M/$s_t$ system, computed using a  shifted $\rho_t$, i.e. a
quantity $\hat \rho_t = \rho_t - \delta$ for appropriate $\delta$.

In a regime where $\rho_t$ is close to (but smaller than) $1$, one
could use one of the popular measures of queueing system performance
as ``cost'' -- or use $\rho_t$ itself as the cost.

\subsection{Other cost functions}
Other situations of practical relevance besides the two considered above
are likely to arise.  Our robust planning methodology,
described below, is flexible and rapid enough
that many alternate models could be
accommodated.  Moreover, we would argue that in the context of the cost
associated with staff shortfall, any reasonable cost function would be,
broadly speaking, increasing as a function of the shortfall.  The approach
described above approximates cost, in the two cases we listed, using piecewise-linear
convex functions, and we postulate that many cost functions of
practical relevance can be successfully approximated this way.

\section{Robust Models}\label{sec:robust}

We now describe our methodology for the robust surge staffing problem.
We will first describe our uncertainty model, then the deployment policies
we consider, and finally the robust optimization model itself.  In what follows, we make the assumption that the total quantity of surge staff
is small enough that their deployment does not affect the evolution of the epidemic
in the SEIR model, in particular the values $\beta_t$.

\subsection{Uncertainty models}\label{subsec:UncertModels}

Let $p_t$ denote the probability of contagion at time $t$ (time measured relative to
the declaration of the epidemic).
The model for uncertainty in $p$ that we consider embodies the notion
of increasing uncertainty in later stages of the epidemic.  We will
fist describe the model and then justify its structure. We assume that
there are four given parameters $p^L, p^U, \hat p^L,$  and $\hat p^U$.  The model behaves as follows:\\

\vspace{.1in}
\noindent There is a time period $\breve t$, unknown to the planner, such that
for $1 \le t \le \breve t$, $p_t$ assumes a constant (but unknown
to the planner) value in $[p^L, p^U]$, and
for $\breve t < t$, $p_t$ assumes a constant (and also unknown
to the planner) value in $[\hat p^L, \hat p^U]$.\\
\vspace{.1in}

\noindent We stress that $\breve t$ is not known to the planner.  We are interested in cases where $\hat p^L \le p^L$ and
$ p^U \le \hat p^U$, that is to say, the period following day $\hat t$
exhibits more uncertainty than the initial period.  We can justify our
model as follows.  In the event of an epidemic, a planner would be able to
obtain some information on the spread of the epidemic in the period leading to the actual declaration of the epidemic, resulting in a (perhaps tight) range
of values $[p^L, p^U]$.  As we will discuss below (Section \ref{subsec:deploymentstrat}) we are interested in surge staff deployment strategies with limited
flexibility, i.e. requiring staff commitment levels that are prearranged.

In this setting,
at the point when the epidemic is declared, the planner would deploy a
staff deployment strategy.  A prudent planner, however, would not simply
accept the range $[p^L, p^U]$ as fixed.  In particular, if at first
the epidemic is mild (small $p^U$) the planner might worry that a change,
such as sudden drop in temperature, could effectively increase $p$ beyond
$p^U$. Note that a change in weather would not
simply change the contact rates, $\lambda;$ it might produce other environmental
changes (such as decreased ventilation); further, it is known that the influenza virus has higher transmissibility in colder and drier conditions \cite{convRaulRabadan, flu_coldweather}.  we model such changes, collectively,
through changes in $p$.  Should the change
take place, a staff surge strategy that consumed most available staff in
the earlier part of the epidemic would be ineffective.

As a means to avoid this overcommitment of resources to early phases of the
epidemic, we can assume that a
second and more virulent regime of the epidemic \textit{could} be manifested
at an unknown later time. We can parameterize this later regime by
assuming a range $[\hat p^L, \hat p^U]$ with, for example, $\hat p^L = p^L$ and $\hat p^U > p^U$.  The exact relationship between the two values
$\hat p^U$ and  $p^U$ is a measure of the conservativeness or risk-aversion of the decision-maker. We will touch upon this issue later.

Of course, exactly the reverse could happen: the epidemic could become
\textit{milder} in the second stage. This could take place as a function of
changes in public behavior due to non-pharmaceutical interventions (see
 \cite{bridges2003transmission, grayson2009efficacySWandHH}) resulting
in a (difficult to accurately predict) decrease in infectivity.  In
that case, a surge strategy that defers most staff deployment to the
second stage of the epidemic could also be ineffective. The challenge is how to properly hedge in view of these extreme situations,
and other intermediate situations that could also arise.

\subsection{Deployment strategies}\label{subsec:deploymentstrat}
We make several assumptions that constrain the feasible deployment patterns.
First, we assume that a surge staff member, when deployed, will be available for up to $\tau$ time periods, provided that he or she does not get infected first - if that happens, this person will be removed from
the system and will not available for deployment in the future.  We also assume that the pool of available
surge staff over the planning horizon is finite.  Finally, there is a maximum
quantity of surge staff that can
be summoned on any given time period.  More detailed models are possible and
easy to incorporate in our optimization framework.

In this paper we will focus on\textit{ offline}, or \textit{fixed} strategies.
More precisely, we assume that a \textit{deployment vector}
is
computed immediately after an epidemic is declared, on the basis of the
available information.  From a formal perspective, the deployment vector
is obtained as follows:
\begin{itemize}
\item [(1)] First, the epidemic must be declared. As indicated in Section
\ref{subsubsec:declare}, this takes place as soon as the percentage of infected individuals
exceeds some (small) threshold.
\item [(2)] Based on data available at that point, estimates are constructed for the
various SEIR parameters.  In particular, the ranges $[p^L, p^U]$ and
$[\hat p^L, \hat p^U]$ for $p$ (see Section \ref{subsec:UncertModels}) are constructed.
\item [(3)] Using as inputs the
population statistics, the cost function, the nominal
SEIR model parameters, and the
uncertainty model for $p$, we compute the deployment vector $h = [h_1, h_2, \ldots, h_{T}]$, where $h_t$ indicates the number of staff procured $t$ time periods {\bf after} the declaration of the epidemic. The parameter $T$ is chosen large enough to encompass one epidemic wave and will be discussed later.
\end{itemize}

\noindent We next address some points implicit in (1)-(3).  First,
recall that in the SEIR framework the formal epidemic will have a starting point
which precedes the time period when an epidemic is actually declared.
The exact
 magnitude of this run-up period is {\em not} known to a decision maker.
In all our notation below (as in (2) above), ``time period 1'' refers to the time period where the epidemic was declared, that is to say, we always
label time periods by the amount elapsed since the declaration of the epidemic.
  When using the SEIR machinery to {\em simulate} an epidemic, of course,
we always proceed as in the formal model, and we merely relabel as ``time = 1''
that period where the declaration condition is first reached (since that
would be the start of the epidemic as far as a planner is concerned).

 Also,
we assume that the epidemic is {\em correctly} declared, that is to say,
the first day in which the criterion in Section \ref{subsubsec:declare} applies
is correctly observed.  In practice, this observation would include
noise (for example, due to individuals infected by a different strain) but we make the
assumption
that the resulting decrease in strategy robustness is small.

Item (2), the estimation of the SEIR model and in particular the
construction of the initial range $[p^L, p^U]$, gives rise to a host of other issues
that are outside the scope of this paper (but are nonetheless important).
We assume that at the time the epidemic is declared there is enough data
(however incomplete, and noisy) to enable the application of $robust$ least
squares methods (see e.g. \cite{elghaoui}) to fit an SEIR model, and that the
model is sufficiently accurate to permit point estimates for all parameters
other than $p$.  Having constructed the initial range $[p^L, p^U]$, the
later range $[\hat p^L, \hat p^U]$ is constructed on the basis of (a) risk
aversion, and (b) environmental considerations.  The role of (a) is clear -for example, we could obtain $\hat p^U$ by adding to $p^U$ a multiple (or fraction)
of the standard deviation of $p$ in the observed data.  As an example for (b), an
epidemic that starts in late-Autumn might possibly become more infectious as
colder weather develops.

Another point concerns the time horizon, $T$, for the SEIR model (see
Section \ref{sec:SEIR}) which to remind the reader is measured from the
point that the epidemic is declared.  It is
assumed that $T$ is large enough to handle an epidemic of interest.  If we assume a fixed (but unknown) $p$, then the milder (i.e., less infective) epidemics
will tend to run longer, but, significantly, will also be less disruptive.
Higher values of $p$ give rise to sharper epidemics in the near term.  From
a technical standpoint, the parameters $p_L$ and $\hat p_L$ can be used to
construct estimates for $T$. In order to be sure that we capture as much of an
epidemic as possible, in this study we  will use values for $T$ up to 150,
which according to our experiments seems more than adequate.

One could consider models of epidemics spanning, for example,
 a one-year period,
with multiple epidemic ``waves''.  If some of these waves are very prolonged
then they will necessarily be mild waves for long periods of time.  As a result,
during such long, mild epidemic periods (1) the social cost will be low,
and (2)
a significant pool of the population will become sick and as a result become
immune to the virus.  Consequently, future epidemic waves will necessarily
be milder and cause decreased social cost (so
surge staff will be less critical).
On the other hand, over a long period we could have  well-separated
strong epidemic waves.  However, we would argue that from the perspective of
surge staff deployment, each such wave should be
handled as a separate event, with its own surge staff deployment strategy.\\

A significant element in our approach is that it produces a fixed deployment
vector $h_1, \ldots, h_T$ which, from a formal standpoint,
will be applied regardless of observed conditions.  Of course,
the strategy
should be interpreted as a template and a planner would
apply small deviations from the planned surge levels, as needed.  In any
case, we have focused on this assumption
for practical reasons.  Having a pre-arranged schedule greatly simplifies the
logistics of deploying possibly large numbers of staff, especially if many
of the staff originate from geographically distant sources.

A broader class of models would allow for on-the-fly revision of a strategy in the midst of an epidemic, which would allow a planner to dynamically react to changes in the epidemic.

However, a significant underlying issue concerns the actual flexibility that
would be possible under a virulent epidemic.   We expect that large,
sudden changes in deployment plans may be difficult to implement. In particular,
if a significant and unplanned \textit{increase} in surge staff is needed, it may prove impossible to
rapidly attain this increase; and by the time the surge staff is available a severe peak of the epidemic may already
have taken place.  See Section \ref{outofsample} for some experimental validation of these
views.

A different type of dynamic strategy would implement many, but small corrections as conditions change.  We will discuss examples of such strategies (and appropriate modifications
to our methodology) in Section \ref{thefuture}; an issue in this context is
the callibration of our statement  ``as conditions change'' in an environment
that will be characterized by noisy, partial and late data.

\subsection{Robust problem}\label{robustproblem}
We can now formally state the problem of computing an optimal robust
pre-planned staff deployment strategy. Let $h$ denote a deployment vector and
$\cH$, the set of allowable deployment vectors. Let $\vec{p} = (p_1, p_2, \ldots, p_T)$ be a vector indicating a value
of $p$ for each time period.  Define
\begin{eqnarray}
V(h \, | \, \vec{p}) & := & \mbox{cost incurred by deployment vector $h$, if the
infection probability equals $p_t$ at time $t$.} \label{eq:costdef}
\end{eqnarray}
\noindent [Here we remind the reader that time is measured relative to the declaration of the epidemic.] Let $\cP$ indicate a set of vectors $\vec{p}$ of interest; our uncertainty set.  Our robust problem can now be formally stated, as follows:

\begin{center}
  \fbox{
    \begin{minipage}{0.5\linewidth}
\begin{center}{\bf Robust Optimization Problem}\end{center}
\vskip -.2in
\begin{eqnarray}
V^* & := & \min_{ h \in \cH} \max_{ \vec{p} \in \cP} ~ V(h \, | \, \vec{p}). \label{eq:V_star_def}
\end{eqnarray}

\end{minipage}
  }
\end{center}

\noindent Problem (\ref{eq:V_star_def}) explicitly embodies the adversarial
nature of our model.  Given a choice of vector $h$, a fictitious adversary
chooses that realization of the problem data (the contagion probabilities) that
maximizes the ensuing cost; the task for a planner is to minimize this
worst-case cost.  The set $\cH$ describes how staff deployment plans are
constrained.  In our implementations and testing, we assume that any deployed staff
person works for a given number $\tau$ of time periods, or until
he or she gets sick; following this period of service this person will
not be available for re-deployment.  We also assume that whenever surge-staff
is called up, there is a one-period lag before actual deployment (this
assumption is easily modified to handle other time lags).

Finally,
we model the impact of the epidemic on surge staff using the SEIR model
with parameters as for the workforce group ($\lambda^2_t$ and $\mu_{E_2}$),
and we assume that the total quantity of available surge staff is small enough
so as to not modify the outcome of the epidemic.

In Appendix \ref{appendix} we will present an efficient procedure,
based on linear programming, for solving problem (\ref{eq:V_star_def}) to
very tight numerical tolerance.

\section{Numerical Experiments}\label{sec:experiments}

In this section we investigate the structure of the optimal robust policies by numerical experiments. First, it is worth raising a point concerning ``approximate'' vs ``exact'' solutions.  Our algorithms
construct numerically optimal solutions to the formal problems they solve;
nevertheless, as argued above, our models approximate real problems.  Again,
we postulate that given the context (i.e. the behavior of social entities)
an approximation is the best we can hope for; moreover we expect that
robust strategies computed for the approximate models will translate into
actual robustness in practice.  This last feature can at least be
experimentally tested through simulation, as we present via examples in this section.

In order to assign numerical values to the various parameters in our models,
we followed the existing literature and consulted with various experts \cite{convRaulRabadan, ConvStephenMorse}).
The numbers in the literature may vary significantly depending on various factors such
as the location where the study took place, the underlying model for which the
parameters were calibrated, and the epidemic wave of study.
For example, published values for the influenza $R_0$,
the number of secondary transmissions caused by an infected individual in a susceptible population, range from $1.68$ to $20$ \cite{mills2004transmissibility}.
Our main sources are listed in Table \ref{tab:RefParamValues}.

\begin{table}[htbp!]\small
  \centering
  \caption{List of references for different parameter values.}
    \begin{tabular}{cccc}
    \addlinespace
    \toprule
    {\bf Parameter} & {\bf Description} & {\bf Value(s)} & {\bf Reference}\\
    \midrule

    \multicolumn{ 4}{l}{{\bf Epidemic related parameters}} \\[3pt]
    $p$   & Probability of contagion & $[0.02, 0.16]$$^a$ & \cite{wallinga2010optimizing}\\ 
          &       & $\{0.07,0.09,0.12\}^b$ & \cite{larson2007infprogressmodel} \\
    $\lambda$ & Contact rate & $47.48$ (per day) & \cite{wallinga2006using}  \\
          &       & $35$ (per day)$^c$ & \cite{larson2007infprogressmodel}\\
    $\mu_E^{-1}$ & Latent period & $1$ day & \cite{fergusonlongini2008layeredcontainment}\\
          &       & $1.9$ days$^d$ & \cite{longini2004contWantiviral}\\
          &       & $1.25$ days & \cite{gardam2007prophylhcw}  \\
          &       & $1.48 \pm 0.48$ days & \cite{gardam2007prophylhcw} and within \\ 
    $\mu_R^{-1}$ & Infectious period & $2$ days & \cite{fergusonlongini2008layeredcontainment}\\ 
          &       & $4.1$ days & \cite{longini2004contWantiviral, lee2007effectiveness}\\
          &       & $4.1$ days$^e$ & \cite{lee2007effectiveness}\\
    $1 - f$ & Mortality rate & $0.02$ & \cite{longini2004contWantiviral}\\
          &       & $0.002$ per day$^f$ & \cite{gardam2007prophylhcw}\\
          &       & $1.2$\%$^g$     & \cite{chowell2011H1N1Mexico}\\
          &       & $0.05$$^g$ & \cite{lee2007effectiveness}\\
    $R_0$ & $R_0$ & $1.68$ to $20$$^h$ & \cite{mills2004transmissibility} and within \\
          &       & $1.73$$^i$ & \cite{wallinga2006using}\\
          &       & $2.5$$^j$ & \cite{lee2007effectiveness}\\
          &Length of pandemic wave & $8 - 12$ weeks & \cite{gardam2007prophylhcw} \\[6pt]

    \multicolumn{ 4}{l}{{\bf Hospital related measures}} \\[3pt]
    $\zeta$ & Patient arrival rate & 2,800 ILI cases/day$^k$ & \cite{lee2007effectiveness} \\ 
    $\mu$ & Service rate & $24 - 30$min/patient$^l$ & \cite{green2004improvingER} \\
    $\rho$ & Occupancy rate & Average of $85\%$$^m$ & \cite{Rutgers2007NJHospitalCapacity}\\ 
    $\delta$ & Increase in demand for health care & up to $50$\% & \cite{NJPandemicPrep}\\[6pt] 

    \multicolumn{ 4}{l}{{\bf Miscellaneous}} \\[3pt]
    & Seasonal threshold   & $[0.06, 0.076]$ weekly$^n$ & \cite{url:CDCFluWeekly}\\
    & Epidemic threshold   & $[0.064, 0.08]$ weekly$^n$ & \cite{url:CDCFluWeekly}\\
    & ILI National Baseline   & $2.4\%$ weekly$^o$ & \cite{url:CDCFluWeekly}\\%
    $\theta$ & Transmission reduction factor  & $29$\% to $37$\%$^p$ & \cite{chowell2011H1N1Mexico}\\
    \bottomrule
    \multicolumn{4}{l}{{\footnotesize $^a$ 1957 wave in the Netherlands}}\\[-1pt]
    \multicolumn{4}{l}{{\footnotesize $^b$ Values of $p$ for different socially active groups}}\\
    \multicolumn{4}{l}{{\footnotesize $^c$ Weighted average of different age groups presented in the paper}}\\
    \multicolumn{4}{l}{{\footnotesize $^d$ 1957 pandemic wave}}\\
    \multicolumn{4}{l}{{\footnotesize $^e$ Untreated symptomatic}}\\
    \multicolumn{4}{l}{{\footnotesize $^f$ Treated symptomatic}}\\
    \multicolumn{4}{l}{{\footnotesize $^g$ Overall ILI case-fatality ratio during 2009 H1N1 pandemic in Mexico}}\\
    \multicolumn{4}{l}{{\footnotesize $^g$ Spanish Flu fatality rate}}\\
    \multicolumn{4}{l}{{\footnotesize $^h$ From multiple studies}}\\
    \multicolumn{4}{l}{{\footnotesize $^i$ 1957 wave}} \\
    \multicolumn{4}{l}{{\footnotesize $^j$ Min 1.5, Max 6}} \\
    \multicolumn{4}{l}{{\footnotesize $^k$ Data for Singapore, 2005. Population then was about 4.35 million people}}\\
    \multicolumn{4}{l}{{\footnotesize $^l$ Emergency Department, New York Presbyterian Hospital}}\\
    \multicolumn{4}{l}{{\footnotesize $^m$ Maintained beds capacity for New Jersey Hospitals, 2005}}\\
    \multicolumn{4}{l}{{\footnotesize $^n$ The threshold is for percentage of deaths caused by pneumonia and influenza.}}\\[-1pt]
    \multicolumn{4}{l}{{\footnotesize It is different for each week throughout the year.}}\\
    \multicolumn{4}{l}{{\footnotesize $^o$  Percentage of influenza-like illness (ILI) patient visits reported through the U.S. Outpatient ILI Surveillance Network.}}\\
    \multicolumn{4}{l}{{\footnotesize $^p$ Estimated effect of social distancing measures during H1N1 epidemic in Mexico, 2009}}
    \end{tabular}
  \label{tab:RefParamValues}
\end{table}

\subsection{Health care setting}\label{subsec:NumExp_Healthcare}

We now present examples modeling a hospital of the size of New York-Presbyterian Hospital in New York City under slightly different scenarios. Staff at this hospital amounted to $19,376$ in 2010 (we have rounded this up to $20,000$) and served a population of
approximately 900,000 people \cite{NYPresbHospital}. We used these numbers as a baseline to obtain different scenarios with different initial conditions.

\noindent In the event of an epidemic, hospitals will be likely to experiment surge demand trailing the incidence curve. We use the queueing-based cost function to capture both shifts in demand and workforce availability while the epidemic ensues. This increase in demand for service is modeled by adding a proportional factor of $I^1_t$ to the average arrival rate of patients when we compute the system's occupation rate, $\rho_t$:

\begin{equation} \label{eq:rho_t_xdem}
    \rho_t = \frac{\bar{\zeta} + \delta \, I^1_t}{s_t \mu}.
\end{equation}


\subsubsection{Example 1}\label{subsubsec:Health_Ex1}
For this example we assume 3,000 surge staff members are available, each of which will be called up for at most one period of service lasting one week. If the member
becomes infectious while on call, then service is terminated with no further future availability.

We used the uncertainty model in Section \ref{subsec:UncertModels} which
allows $p_t$ to change once. The change is further constrained to take place
within a fixed range of time periods.

The social contact model assumed the nonhomogeneous-mixing rule described in section \ref{subsubsec:NonHomMix_SocDist}. Additionally, we assumed that social
contact
rates were dampened by a factor of $30$\% while the epidemic is declared; that is, during the days in which the growth of infectives is higher than a given threshold. Additionally, we assumed an epidemic is declared only once. In other words, once the rate of infectives slow down below the epidemic threshold, the
epidemic is considered to be terminated.

The parameter values for this example are summarized in Table \ref{tab:ExQParamValues}; the range of days where $p$ is allowed to change are measured from
the \textit{start} of the epidemic (rather than the day the epidemic is
declared).

\begin{table}[htbp]
  \centering
  \caption{Parameter values for example 1}
    \begin{tabular}{ccc}
    \addlinespace
    \toprule
    {\bf Parameter} & {\bf Description} & {\bf Value(s)} \\
    \midrule
    \multicolumn{ 3}{l}{{\bf Epidemic related parameters}} \\
    $P$   & Uncertainty set & $[0.01, 0.012] \times [0.0125, 0.0135]$ \\
          & Possible days of change & $ \{ 140, \dots, 160\}$ \\
    $(\Lambda_1, \Lambda_2)$ & Contact rates & $(30, 35)$ per day \\
    $\mu_E^{-1}$ & Latent period & $1.9$ days \\
    $\mu_R^{-1}$ & Infectious period & $4.1$ days \\
    $1 - f$ & Mortality rate & $0$ \\
    $R_0$ & $R_0$ & $[1.24,1.48]$\\ 
          & Deployment threshold & $2.4\%$ weekly \\[3pt]

    \multicolumn{ 3}{l}{{\bf Population parameters}} \\
    $N_1$ & General population size & 900,000 \\
    $N_2$ & High risk population size & 20,000 \\
    $[ I_1, I_2]$ & Initial infectives & [5,0] \\[3pt]

    \multicolumn{ 3}{l}{{\bf System Utilization (Queueing setting)}} \\
    $\bar \zeta$ & Patient arrival rate & 500 per day \\
    $\mu$ & Service rate per person & $30$ patients per day \\
    $\rho_0$ & Initial occupation rate & 0.875 \\
    $\delta$      & Daily increase in demand for health care & $0.07$\% \\[3pt]

    \multicolumn{ 3}{l}{{\bf Deployment parameters}} \\
                & Total number of available volunteers & $3,000$ \\
     $\tau$     & Length of stay (w/o sickness) & 7 \\
                & Deployment lag & 1 \\
    \bottomrule
    \end{tabular}
  \label{tab:ExQParamValues}
\end{table}

In what follows, the staff deployment plan computed by our algorithm will
be called the \textit{Robust Policy}. To gauge its usefulness, we study its
performance under different scenarios, each of which is defined by a tuple $(p_1, p_2, d)$ representing the initial probability of contagion ($p_1$), the time period where it changes ($d$)
and the value to which it changes ($p_2$).

A natural alternative to the Robust Policy is that which best responds to what could be construed as a worst-case scenario -
the scenario which yields the highest cost when no contingency plan is implemented. Such a strategy can be obtained by solving a linear program as described in (\ref{eq:SEIR_Vhp}) once the tuple that implies the evolution of $p_t$ in this scenario is identified. In this example such tuple is given by $(0.01092,0.0135,140);$ we will term it the \textit{No-Action-Max-Cost tuple},  and
the deployment strategy which achieves minimum cost under this tuple the \textit{Na\"ive-worst-case Policy}. We compare the Na\"ive-worst-case and Robust Policies under the following scenarios:
\begin{enumerate}
    \item No-Action-Max-Cost tuple $(0.01092,0.0135,140)$;
    \item The tuple that achieves highest cost if the \textit{Na\"ive-worst-case Policy} is implemented; and
    \item The tuple that achieves highest cost if the \textit{Robust Policy} is implemented.
\end{enumerate}

\begin{table}[htbp]
  \centering
  \small
  \caption{Results. Example 1.}
    \begin{tabular}{crccc}
    \addlinespace
    \toprule
    Scenario / Strategy / tuple &       & {\bf No Intervention} & {\bf Robust Policy} & {\bf Na\"ive-worst-case Policy} \\
    \midrule
    1     & Cost  & 4.5812 & 0.0495 & 0.0000 \\
    No intervention & Maximum $\rho$ &1.0481  & 1.0018 & 0.9999 \\
    (0.01092,0.0135,140) & \# days $\rho \geq$ 1 & 28    & 8     & 0 \\
         &   &  &  &  \\
    2     & Cost  & 1.4300 & 0.0500 & 0.7100 \\
    Na\"ive-worst-case Policy & Maximum $\rho$ & 1.0210 & 1.0020 & 1.0185 \\
    (0.01172, 0.0135, 140) & \# days $\rho \geq$ 1 & 20    & 8     & 13 \\
         &   &  &  &  \\
    3     & Cost  & 1.6938 & 0.0521 & 0.6862 \\
    Robust Policy & Maximum $\rho$ & 1.0236 & 1.0027 & 1.0170 \\
    (0.01168, 0.0135, 140) & \# days $\rho \geq$ 1 & 21    & 7     & 12 \\
    \bottomrule
    \end{tabular}
  \label{tab:Ex1QResults}
\end{table}

Results are summarized in Table \ref{tab:Ex1QResults}. In Scenario 1 the Robust Policy presents a cost improvement of
almost $99\%$ over the no-intervention policy. While it is not completely able to prevent $\rho$ from exceeding $1$,
it reduces the number of critical days from 28 to 8, and the maximum $\rho$ from $1.05$ to below $1.002$.
Indeed, the undesirable instability in which a system is put through when $\rho > 1$ is heavily penalized with the use of the exponential function in our objective function (see \ref{subsec:queue}). On the other hand, under Scenario 1 the Na\"ive-worst-case Policy reduces the cost to $0$; that is, there is no day in which $\rho$ is greater than $1$.  Without doubt, it represents a better solution than the Robust Policy \textit{for this particular case}. This result is expected; the solution derived
from that single linear program is \textit{specialized} in this scenario, while the Robust Policy,
having to hedge against other possible bad cases, does not perform as well in this instance.


\begin{figure}[h!]
\centering
\subfloat[\textbf{Scenario 1}:Effectiveness of Robust Policy]{
\includegraphics[width=0.45\textwidth]{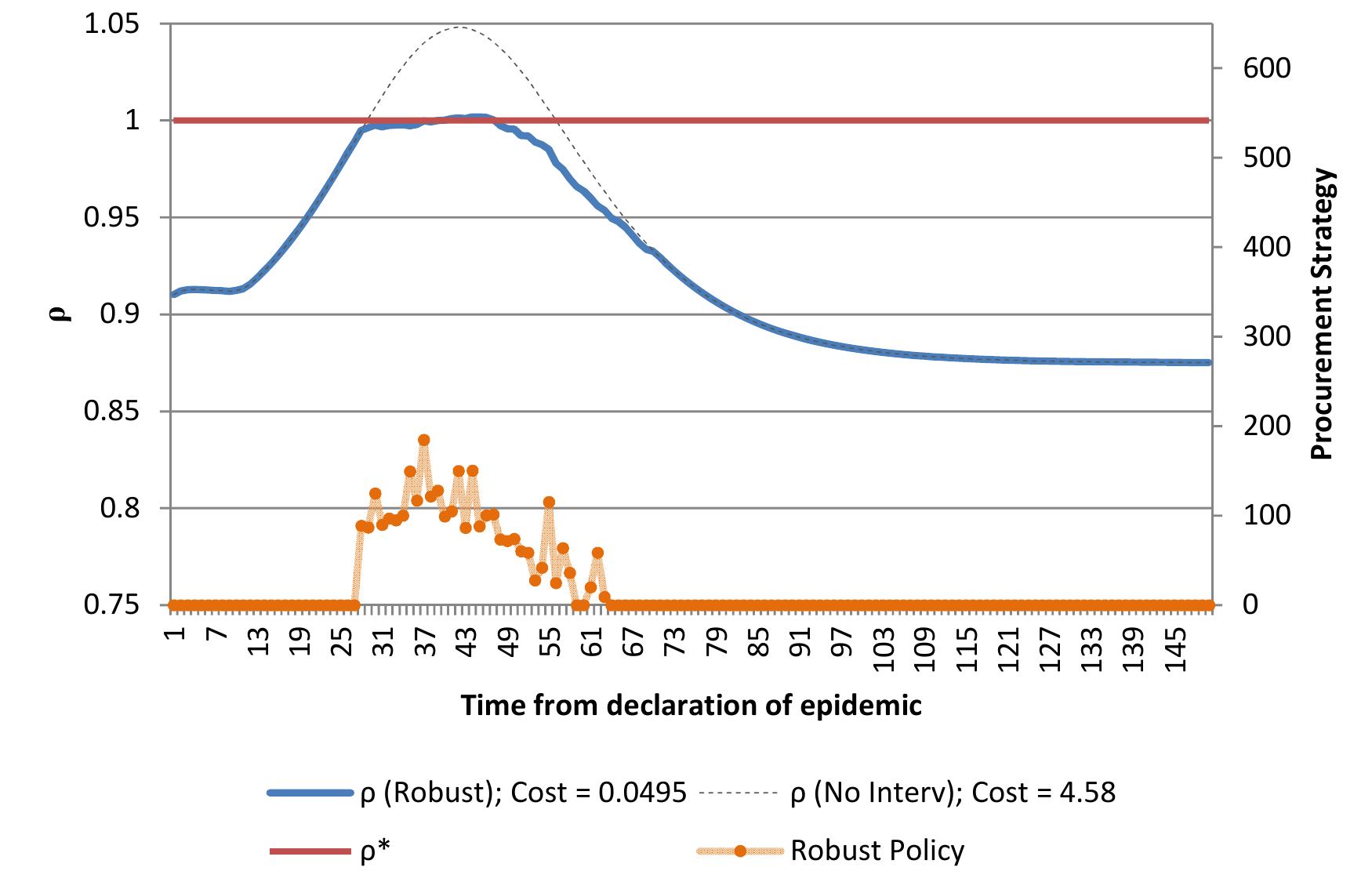}
\label{chart:QR_wp00}
}
\subfloat[\textbf{Scenario 1}:Effectiveness of the Na\"ive-worst-case Policy]{
\includegraphics[width=0.45\textwidth]{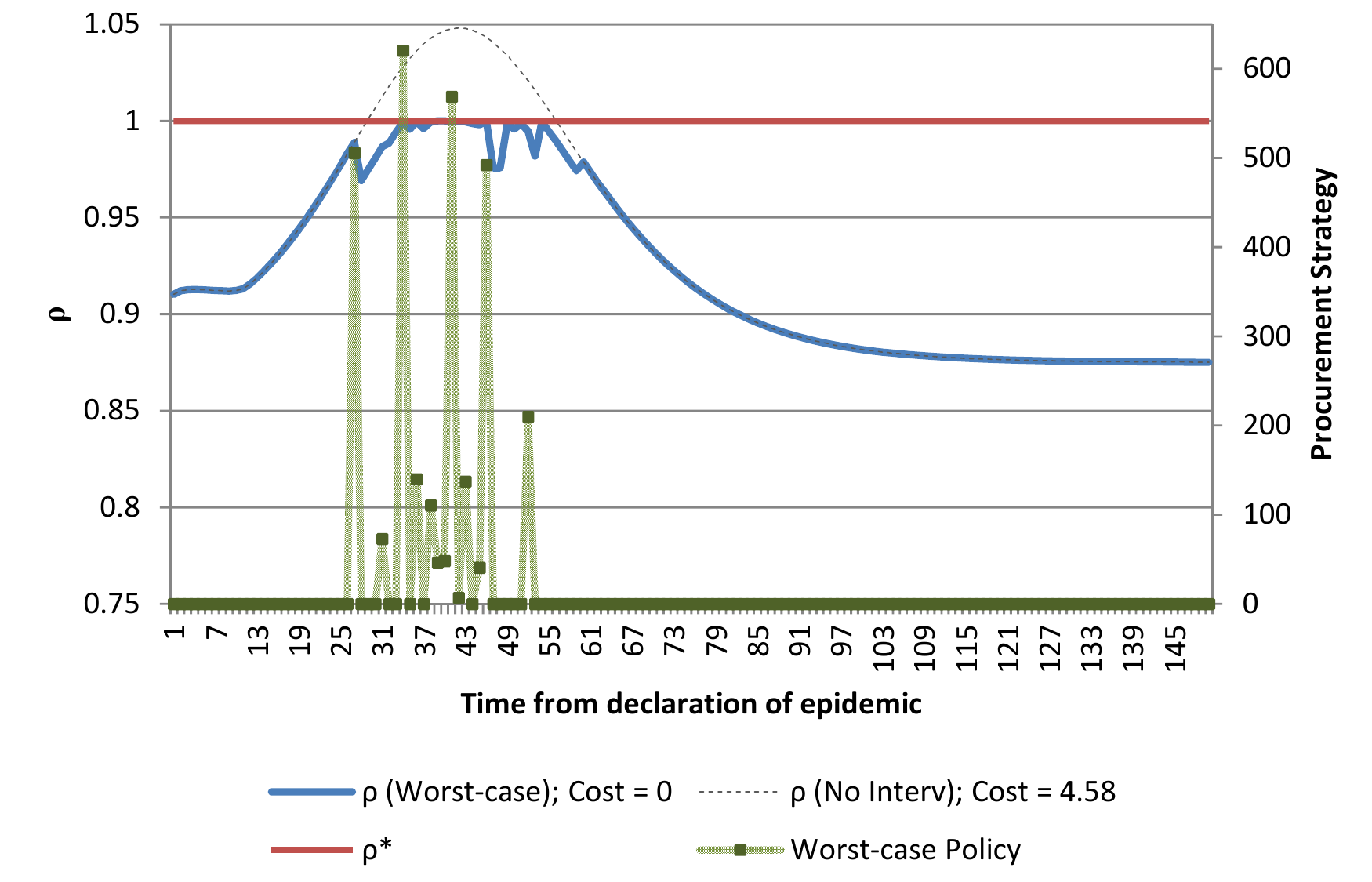}
\label{chart:QLP_wp00}
}
\vspace{10pt}
\subfloat[\textbf{Scenario 2}:Effectiveness of Robust Policy]{
\includegraphics[width=0.45\textwidth]{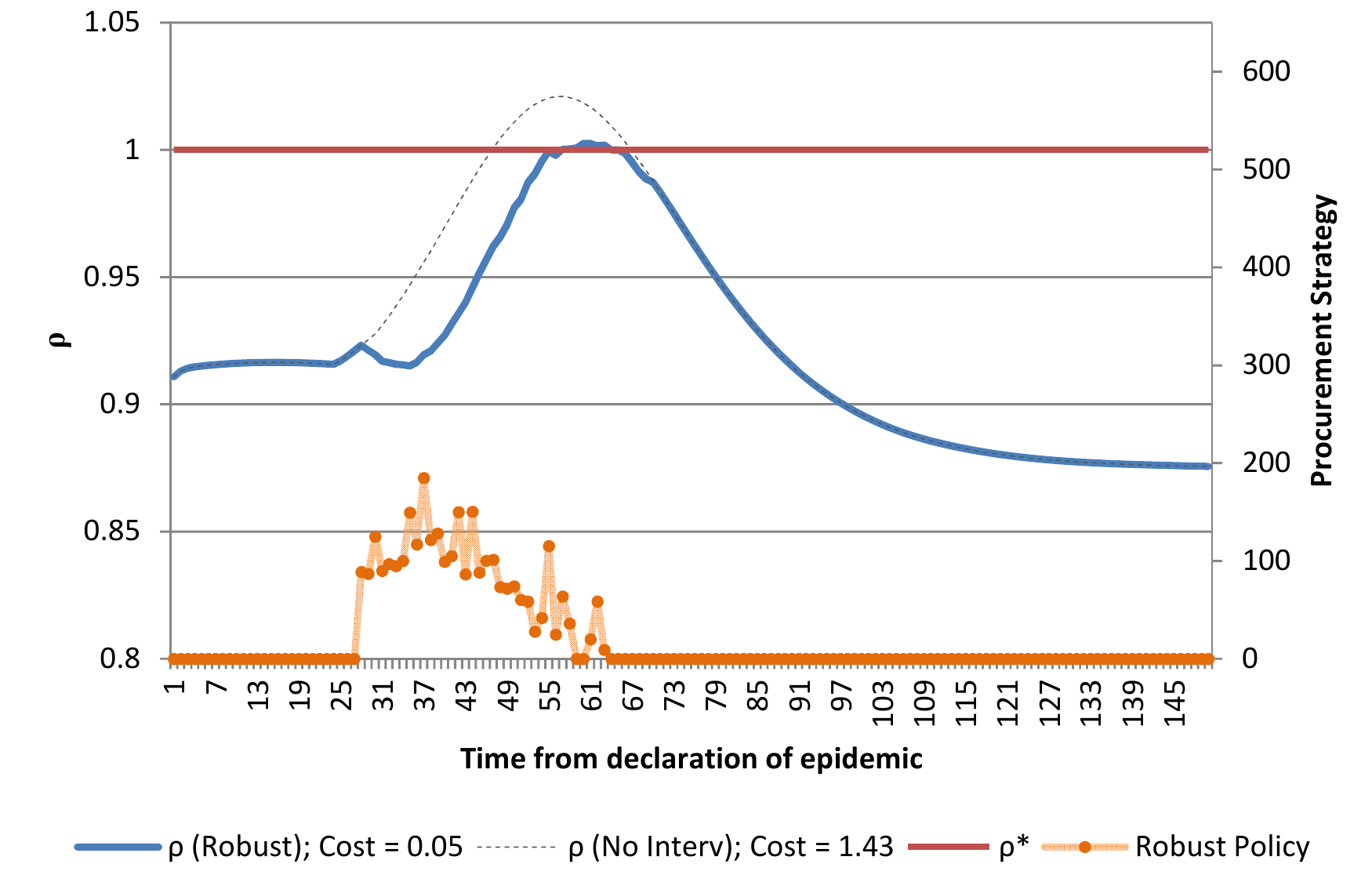}
\label{chart:QR_wp01}
}
\subfloat[\textbf{Scenario 2}:Effectiveness of Na\"ive-worst-case Policy]{
\includegraphics[width=0.45\textwidth]{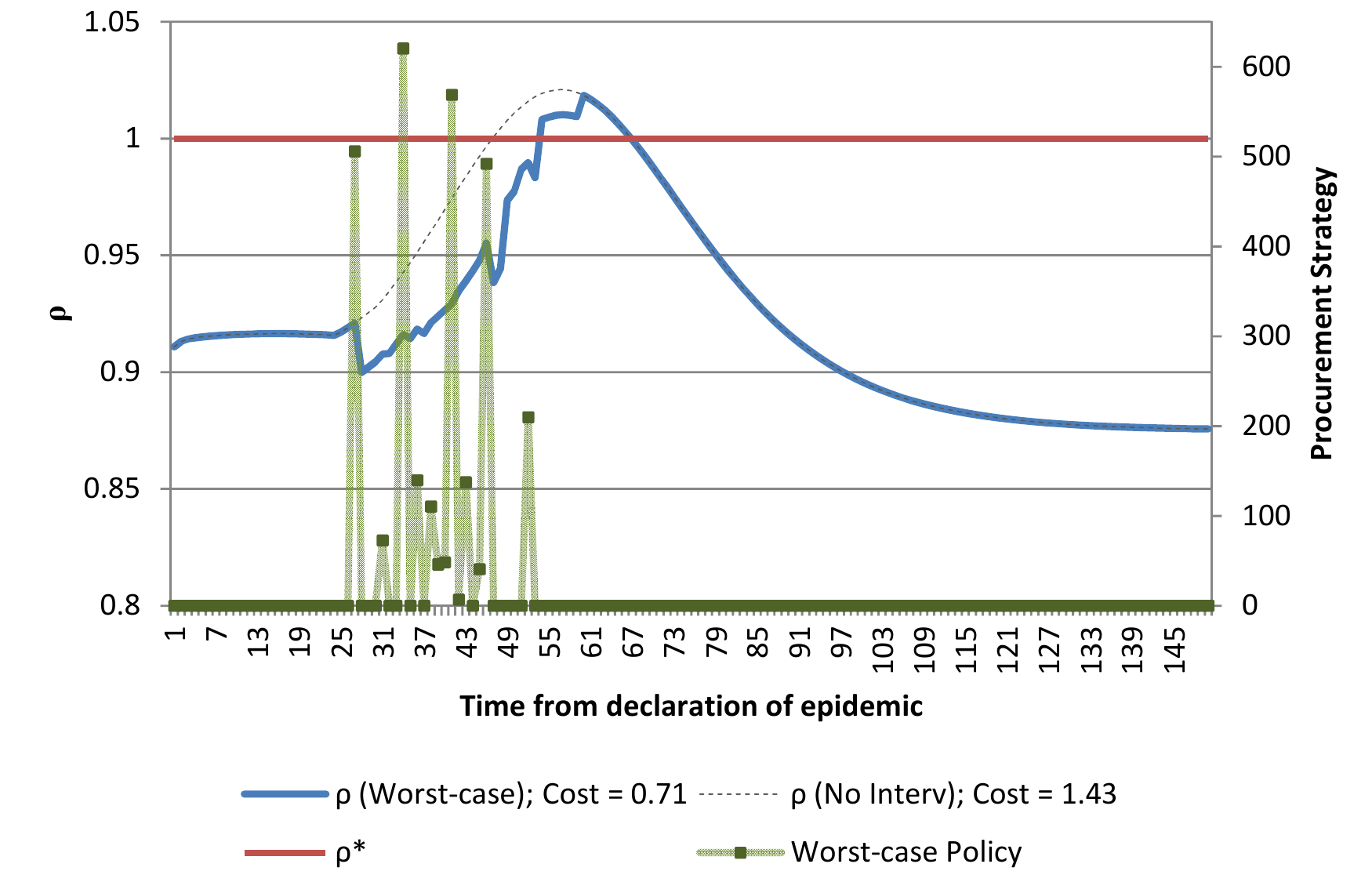}
\label{chart:QLP_wp01}
}
\caption{Reduction in $\rho$ obtained by the Robust and Na\"ive-worst-case deployment strategies under two scenarios: 1) Worst-case scenario given that no policy is implemented, and 2) Worst-case scenario given that Na\"ive-worst-case Policy is implemented.}
\label{chart:Example1}
\end{figure}

The situation is considerably different, however, when we look at Scenario 2, which achieves
the highest cost when the Na\"ive-worst-case Policy is being implemented. In this case, the Robust Policy performs much better
than the Na\"ive-worst-case Policy: the former reduces the cost of this scenario by $96.5\%$, while the latter does so by only $50\%$; the Robust Policy reduces the number of critical days from $20$ to $8$ days, while this number jumps to $13$ for the Na\"ive-worst-case Policy. The maximum value of $\rho$ is $1.002$ for Robust Policy and $1.018$ for Worst-Scenario, compared to $1.021$ with no intervention.

This fact is further illustrated in Figure \ref{chart:Example1}, where we compare the
behavior of $\rho_t$ with and without interventions for Scenarios 1 and 2 together with the deployment patterns. We note that the Na\"ive-worst-case Policy deploys staff around in large,
sharp bursts between days $27$ and $52$. The Robust Policy deployment, on the other
hand, trails $\rho$ in Scenario 1 in a smoother way by calling in personnel in
smaller batches, and also over a longer period of time: $36$ days.
This shows the importance of Scenario 1 as the costliest case and, at the same time,
 how the strategy hedges against scenarios where the epidemic has its peak later (than under
the No-Action-Max-Cost tuple), even if not so aggressively. In particular,
the last wave of staff deployed by the Robust Policy (on days $54, 58,$ and $63$)
could seem somehow wasteful under Scenario 1; however their utility is displayed under Scenario 2, where the Robust Policy tackles a delayed, weaker $\rho$ peak much more effectively
than the Na\"ive-worst-case Policy. Because an unstable system builds up a queue exponentially fast, even this small changes in $\rho$ could translate in significant difference in the load experienced by the system.

A similar situation is presented in Scenario 3, as illustrated in Table \ref{tab:Ex1QResults}. The Na\"ive-worst-case Policy is not as
effective as handling what would have been a milder epidemic, while the Robust Policy performs consistently well.
In fact, this is the worst scenario for the Robust Policy; nevertheless it manages to significantly outperform the Na\"ive-worst-case Policy here
as well.  Indeed, a feature of the Robust Policy is its near-uniform behavior under all scenarios; this epitomizes the use
of the term 'robust'.

\subsubsection{Example 2}\label{subsubsec:Health_Ex2}

To further illustrate the structure of the Robust Policy under different scenarios, we now consider the previous example
with some modifications.  Here, the uncertainty set is $p$ lies in the interval $[0.01, 0.0125]$ throughout,
but can change values once on any day in the range $\{ 100, \dots, 115\}.$ $R_0$ is now between 1.24 and 1.54. The total number of available volunteers is $2,000$ and the contact rates are not dampened while the epidemic is being declared.
All other data remains as in the first example.

Following the same logic as in the previous example, we present the Robust Policy and compare it to the Na\"ive-worst-case Policy
(obtained as before by solving the linear program associated with the worst-case tuple of our uncertainty set). We present results for the same three scenarios: the No-Action-Max-Cost tuple, the costliest tuple for the Na\"ive-worst-case Policy, and the costliest tuple for the Robust Policy, respectively. This fact hints on the nonconvexity of the problem, a topic we will touch upon later in this section.

\begin{table}[tbp]
  \centering
  \small
  \caption{Results. Example 2.}
    \begin{tabular}{crccc}
    \addlinespace
    \toprule
    Scenario / Intervention / Tuple &       & {\bf No Intervention} & {\bf Robust Policy} & {\bf Na\"ive-worst-case Policy} \\
    \midrule
    1     & Cost  & 3.8332 & 0.0280 & 0.0066 \\
    No intervention & Maximum $\rho$ & 1.0410 & 1.0012 & 1.0003 \\
    (0.0125, 0.0125, 100) & \# days $\rho \geq$ 1 & 27    & 10    & 6 \\
    \\
    2     & Cost  & 3.5906 & 0.0280 & 0.0542 \\
    Na\"ive-worst-case Policy & Maximum $\rho$ & 1.0393 & 1.0019 & 1.0025 \\
    (0.01235, 0.0125, 112) & \# days $\rho \geq$ 1 & 26    & 7     & 8 \\
    \\
    3     & Cost  & 3.7320 & 0.0295 & 0.0379 \\
    Robust Policy & Maximum $\rho$ & 1.0403 & 1.0010 & 1.0017 \\
    (0.01015, 0.0125, 101) & \# days $\rho \geq$ 1 & 26    & 10    & 10 \\
    \bottomrule
    \end{tabular}
    \label{tab:Ex2QResults}
\end{table}

\begin{figure}[tbp]
\centering
\includegraphics[width=0.5\textwidth]{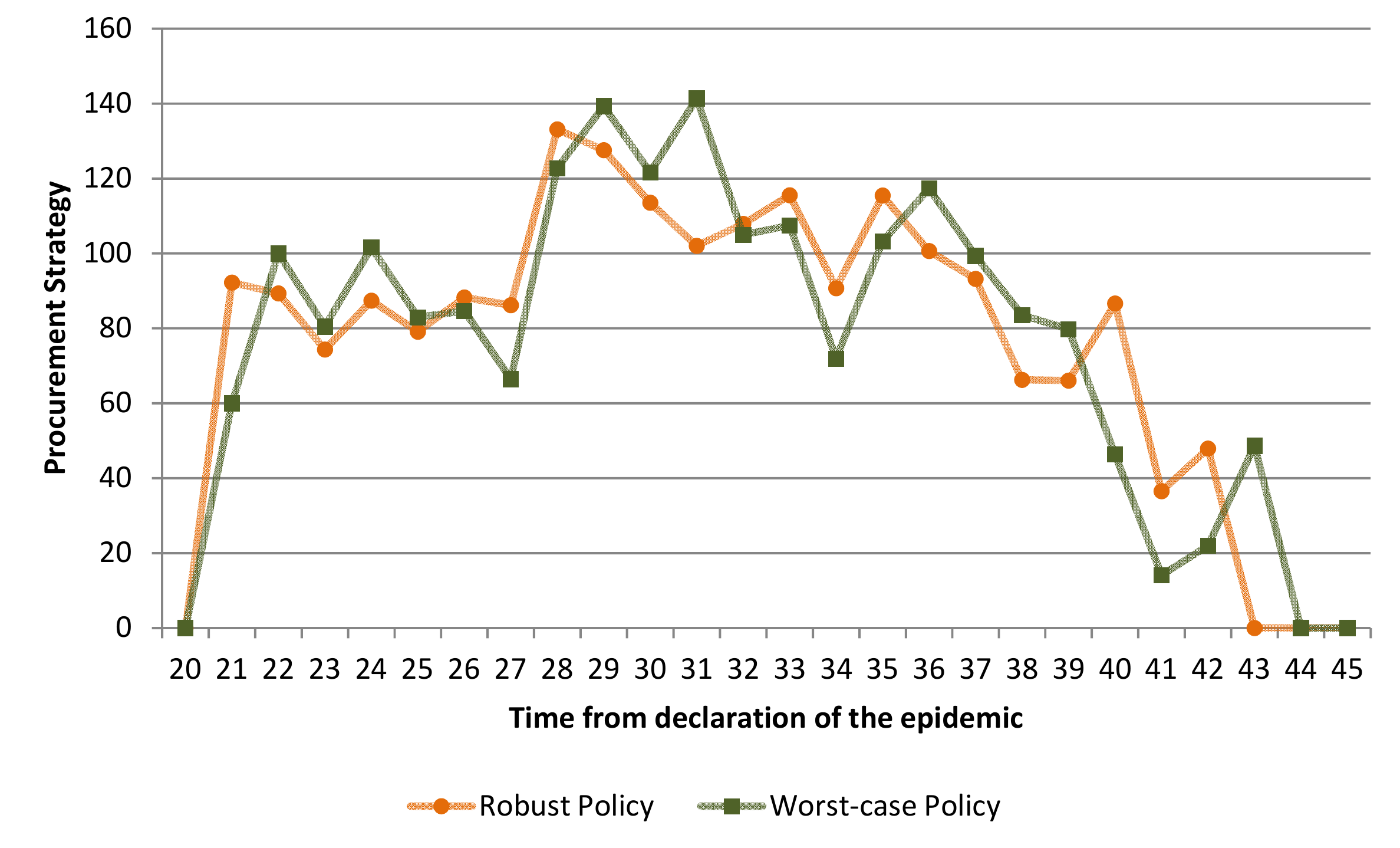}
\caption{Example 2. Deployment Strategies}
\label{chart:Ex2Q_PS}
\end{figure}

Table \ref{tab:Ex2QResults} displays the results. The Robust Policy is worse than the Na\"ive-worst-case Policy in the first instance, while keeping a fairly consistent and low cost for the other scenarios; where the Na\"ive-worst-case Policy does not perform as well.
What is worth noting is that, unlike in the previous example, the structure of the deployment strategies is much more similar (see Figure \ref{chart:Ex2Q_PS}); nevertheless, the small differences make the Robust Policy more resilient to other scenarios.
A second observation regards the structure of the worst-case tuples evinced by the policies.
In this example the No-Action-Max-Cost tuple corresponds to having an epidemic with a \textit{fixed} probability of contagion
$p = 0.0125,$ the highest in the interval of the uncertainty set.
However, the costliest tuples for the Robust and Na\"ive-worst-case policies are quite different.

\subsubsection{Cost-benefit analysis}\label{subsubsec:cost-benefit}

\begin{figure}[htbp]
\centering
\subfloat[Robust Policies.]{
\includegraphics[width=0.45\textwidth]{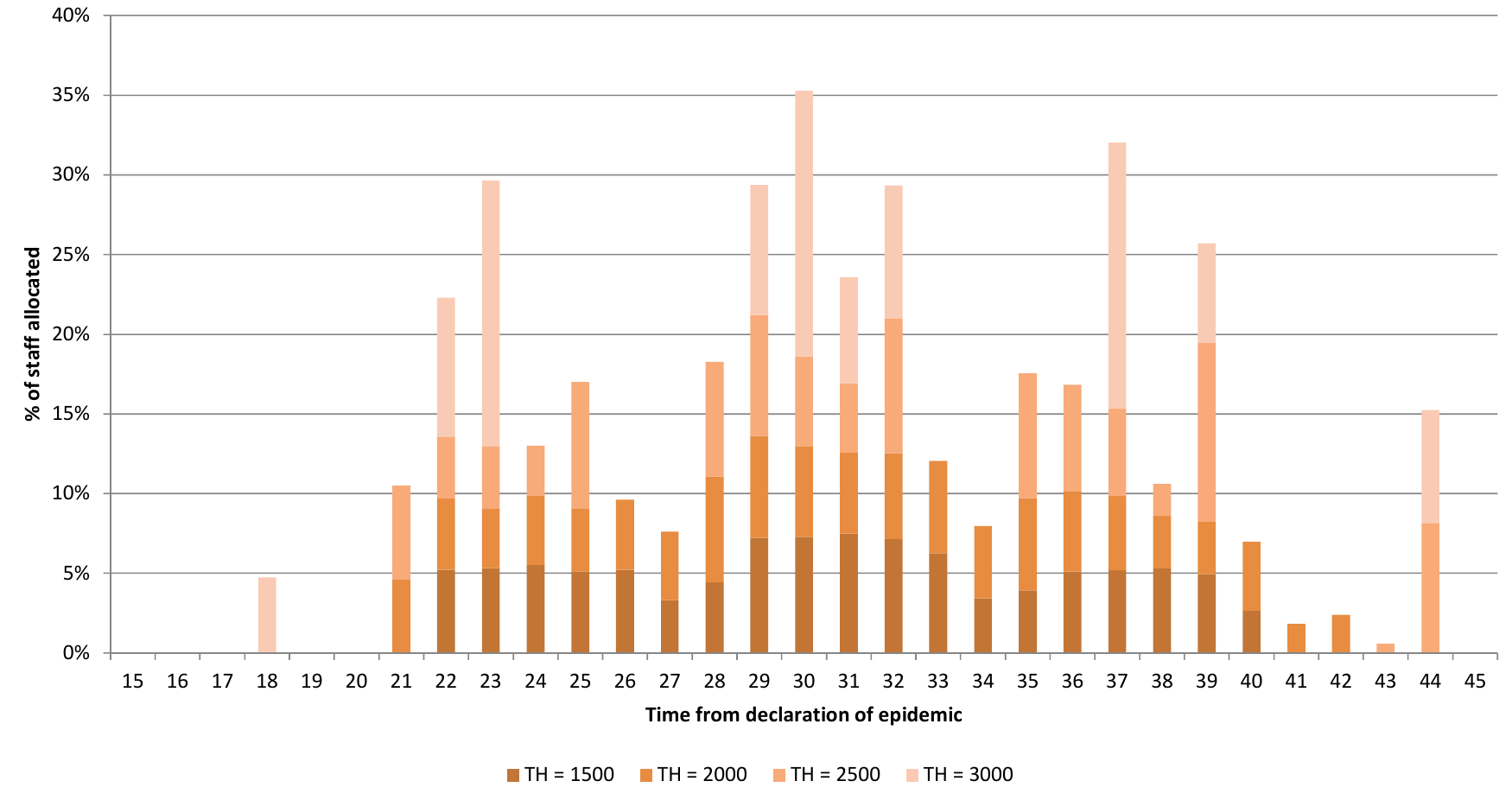}
\label{chart:Ex2CostBenefit_Rob}
}
\subfloat[Na\"ive-worst-case Policies.]{
\includegraphics[width=0.45\textwidth]{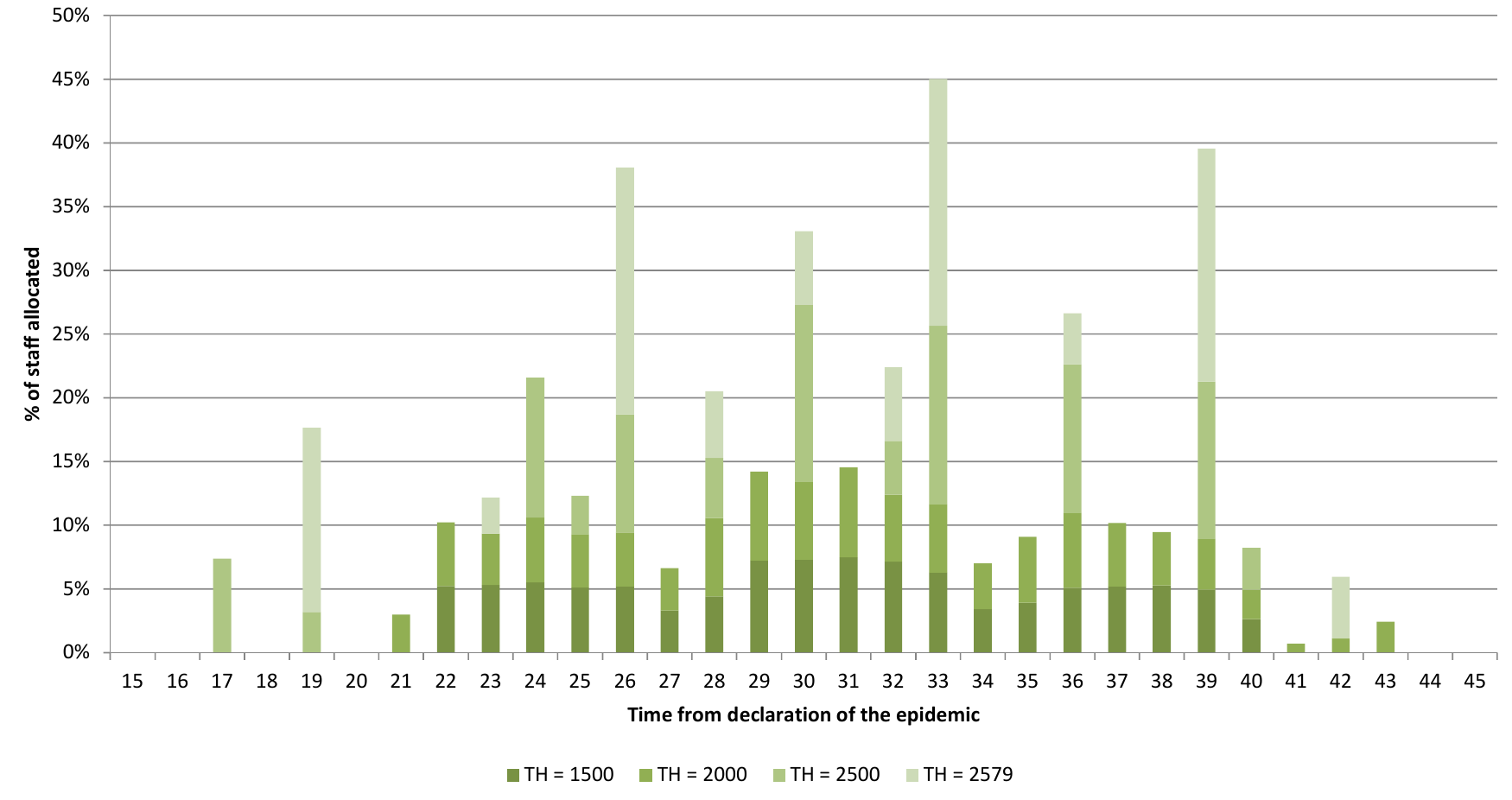}
\label{chart:Ex2CostBenefit_LP}
}
\caption{Distribution of the deployment strategies for different quantities of available surge staff.}
\label{chart:Ex2CostBenefit}
\end{figure}

We now present a cost-benefit analysis on the number of surge staff available. We took this example as a base case and recomputed the optimal Robust and the Na\"ive-worst-case policies assuming $1,000$, $1,500$, $2,500$, and $3,000$ available surge staff.
Table \ref{tab:Ex2_CostBen} displays the change in the worst-case cost given that the corresponding policy is being implemented and its corresponding cost if no contingency plan had been put in place. Similarly, the change in the maximum value of $\rho$ and the number of days that $\rho$ is at or above $1$ are also presented. The last row of the table shows the ratio between the change in the worst-case cost and the change in the number of available surge staff. For example, when there are $1,500$ volunteers available, there is a change in the cost of $0.156\%$ for each of the $500$ that are now not at hand. At the same time, the benefit obtained per additional worker available decreases for both the Robust and the Na\"ive-worst-Case Policy, but more so for the latter. This is again because the Robust Policy allocates the additional staff to cover up for different bad cases, while the Na\"ive-worst-case Policy has only one scenario to focus on. This is depicted in Figure \ref{chart:Ex2CostBenefit}: for the cases in which the Robust Policy has more surge staff available, more people are hired towards the end of the planning horizon, while the Na\"ive-worst-case Policies look to tackle the system congestion early after the epidemic has been declared.

\begin{table}[htbp]
  \centering
  \footnotesize
  \caption{Cost-Benefit Analysis for different quantities of available surge staff.}
  \resizebox{\textwidth}{!}{
    \begin{tabular}{r|ccccc|ccccc}
    \addlinespace
    \toprule
          & \multicolumn{ 5}{c}{{\bf Robust Policy}} & \multicolumn{ 5}{|c}{{\bf Na\"ive-worst-case Policy}} \\
    \midrule
    Total staff deployed & 1,000 & 1,500 & 2,000 & 2,500 & 3,000 & 1,000 & 1,500 & 2,000 & 2,500 & 2,579 \\[3pt]
    Worst-case cost & 1.756 & 0.8183 & 0.0295 & 0     & 0     & 1.7535 & 0.8114 & 0.0542 & 0.0190 & 0.0156 \\[3pt]
    Worst-case cost (No Int) & 3.8332 & 3.8332 & 3.732 & 3.8332 & 3.8332 & 3.8332 & 3.8332 & 3.6077 & 3.6997 & 3.5906 \\[3pt]
    Maximum $\rho$ & 1.02  & 1.02  & 1.001 & 0.999 & 0.998 & 1.02  & 1.01  & 1.002 & 1.002 & 1.002 \\[3pt]
    Maximum $\rho$  (No Int) & 1.04  & 1.04  & 1.04  & 1.04  & 1.04  & 1.04  & 1.04  & 1.04  & 1.04  & 1.04 \\[3pt]
    \# days  $\rho \geq$ 1 & 27    & 27    & 10    & 0     & 0     & 27    & 27    & 8     & 3     & 2 \\[3pt]
    \# days  $\rho \geq$ 1 (No Int) & 27    & 27    & 26    & 27    & 27    & 27    & 27    & 26    & 26    & 26 \\[3pt]
    $\frac{\Delta \mbox{Cost}}{\Delta \mbox{Surge Staff}}$ & 0.345\% & 0.158\% & -     & -0.006\% & -0.003\% & 0.345\% & 0.156\% & -     & -0.002\% & -0.002\% \\[3pt]
    \bottomrule
    \end{tabular}}
  \label{tab:Ex2_CostBen}
\end{table}

\subsubsection{Nonconvexity}\label{subsubsec:cost-nonconv}

As the examples above make clear, the systems we study can be
strongly dependent on the choice of $p$, with clear nonlinearities and nonconvexities.  This justifies the use of robustness in the design of response strategies.

\begin{figure}[tbp]
\centering
\subfloat[Total cost as a function of $p$ under no intervention.]{
\includegraphics[width=0.45\textwidth]{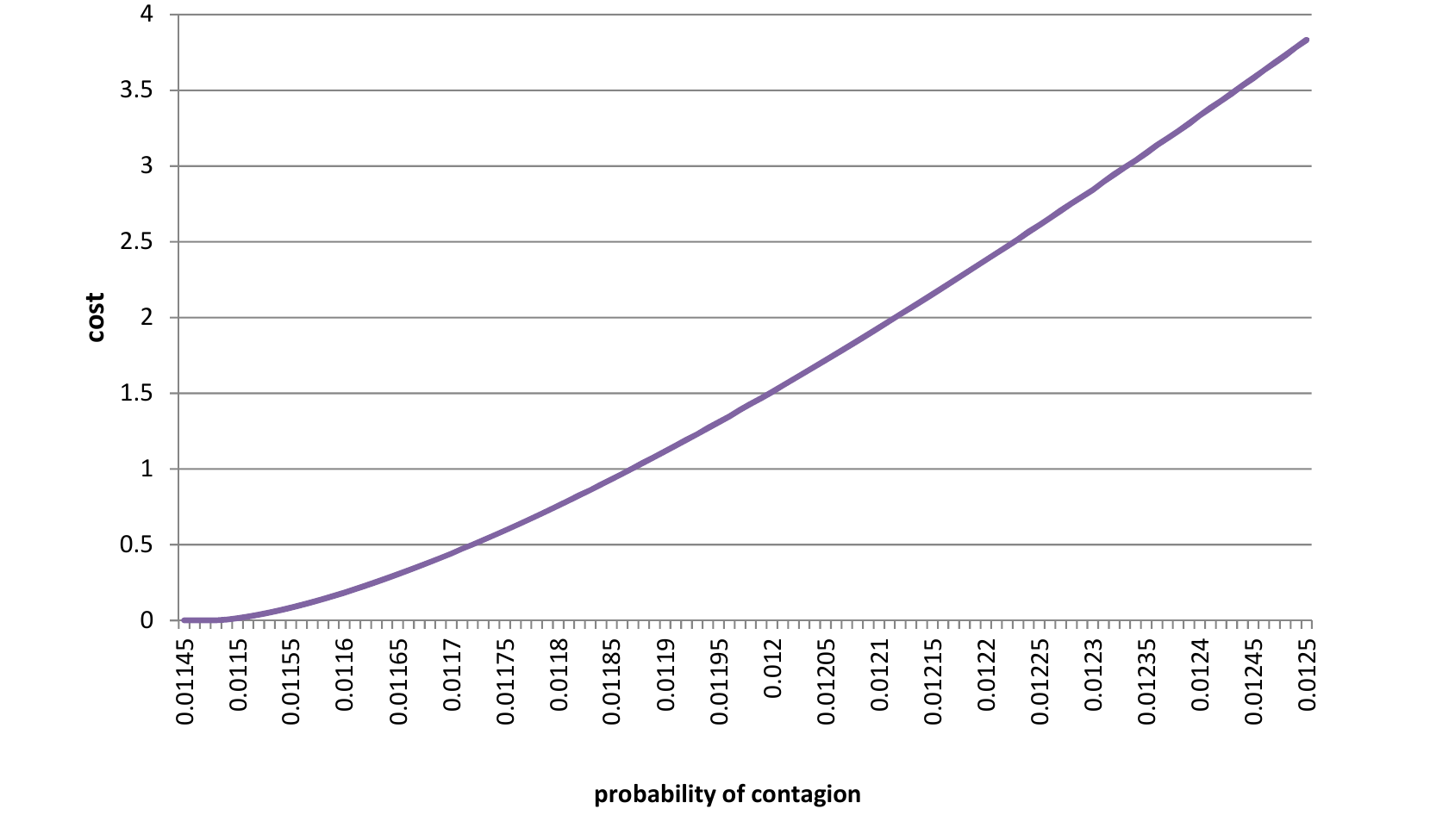}
\label{chart:C_p_NoInt}
}
\subfloat[Total cost as a function of $p$ under 3 different strategies.]{
\includegraphics[width=0.45\textwidth]{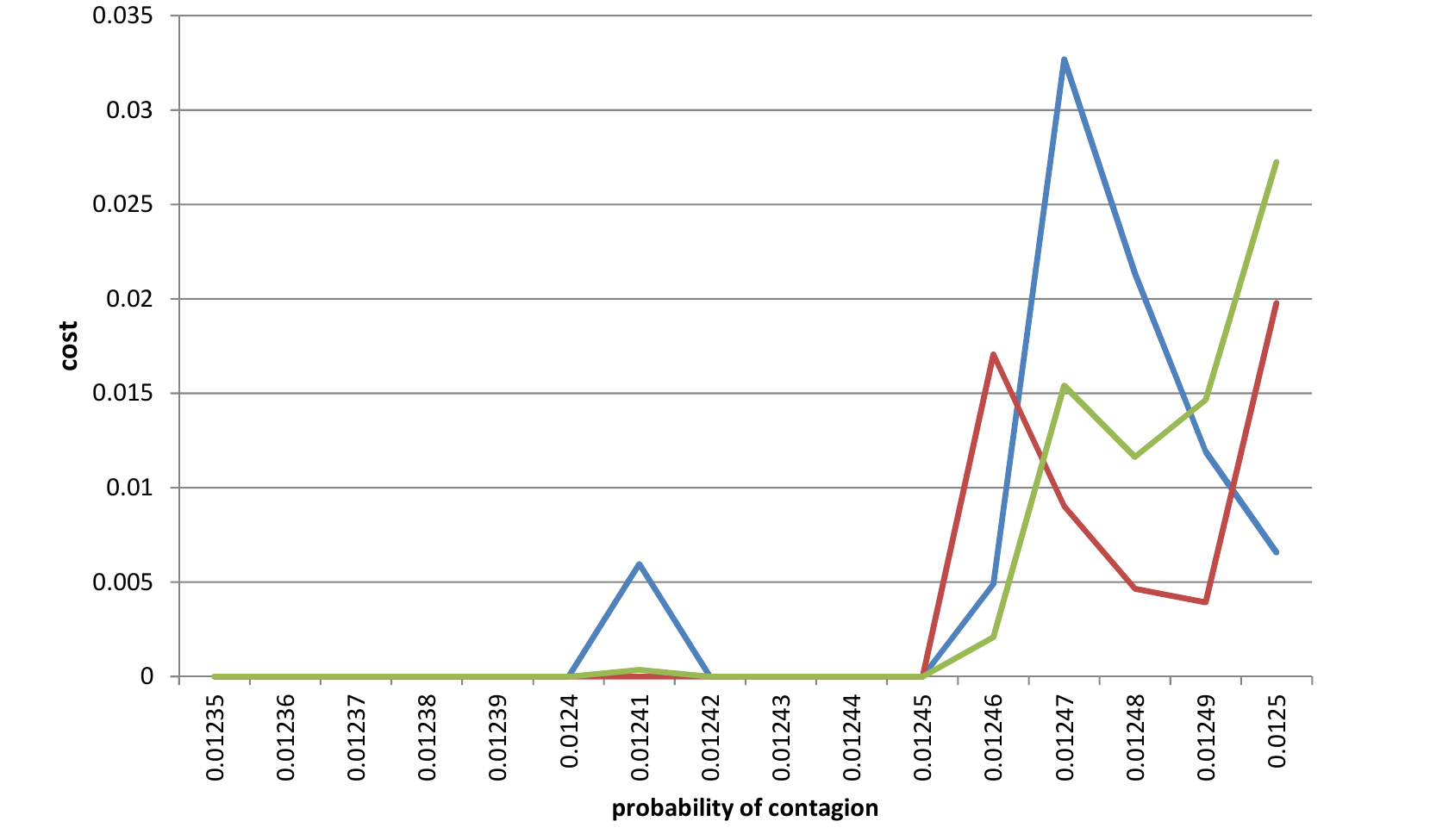}
\label{chart:C_p_3H}
}
\caption{Total cost of an epidemic as a function of $p$ under a)No intervention, and b) Three different deployment policies.}
\label{chart:Ex2_C_p_}
\end{figure}

To further study the dependence on $p$, we considered the cost of a set of strategies when $p$ is constant throughout the epidemic.
Figure \ref{chart:C_p_NoInt} plots the objective function as a function of $p$ in the interval $[0.0115, 0.0125]$ when there is no intervention;  clearly cost is monotone increasing in $p$. However, Figure \ref{chart:C_p_3H} shows similar plots for
three strategies computed by our algorithms\footnote{These strategies were obtained from intermediate iterations of our Robust algorithms. These are fully described in Appendix \ref{appendix}.}.  These
curves are certainly not monotonic and, not surprisingly, particular values of $p$ are able to exploit relative ``gaps'' in deployment.

\subsubsection{Out-of-sample tests}\label{outofsample}
The set of experiments we present below amount to out-of-sample testing.  We build on Example 1 (Section \ref{subsubsec:Health_Ex1}), in which the initial uncertainty interval for $p$ was $[0.01, 0.012]$, and the
second interval was $[0.0125, 0.0135]$; $p$ was allowed to
change between days $140$ and $160$.  In the tests reported in Tables
\ref{tab:late_change1} and \ref{tab:late_change2}, $p$ is allowed to change
to a value $larger$ than $0.0135$ $and$ allowed to change later than in day $160$.  Thus, the experiments go out-of-sample in two ways.

\begin{table}[htbp]
  \centering
  \caption{Effects of a late change in $p$ on Scenario 3}
  \vspace{-1pt}
    \begin{tabular}{r|c|cccc|cccc}
    \multicolumn{ 10}{c}{{ Worst-case scenario given that the Robust Policy is implemented}} \\[3pt]
    \toprule
    {\bf Scenarios} & {\bf Original } & \multicolumn{ 4}{c}{{\bf Hypothetical 1}} & \multicolumn{ 4}{|c}{{\bf Hypothetical 2}} \\
    \midrule
    $p_1$   & 0.01168 & \multicolumn{ 4}{c}{0.01168} & \multicolumn{ 4}{|c}{0.01168} \\
    $p_2$  & 0.0135 & \multicolumn{ 4}{c}{0.014} & \multicolumn{ 4}{|c}{0.015} \\
    day epidemic declared & 113   & \multicolumn{ 4}{c}{113} & \multicolumn{ 4}{|c}{113} \\
    day deployment starts & 140   & \multicolumn{ 4}{c}{140} & \multicolumn{ 4}{|c}{140} \\
    day change $p$ & 140   & 150   & 155   & 160 & 165 -   & 150   & 155   & 160  & 165  \\
    cost Robust Policy & 0.0508 & 0.3268 & 0.0729 & 0  & 0  & 2.1737 & 1.4068 & 0.6282 & 0.0294 \\
    cost No Intervention & 4.58  & 1.6087 & 0.7619 & 0.087 &  0 & 4.1327 & 2.6688 & 1.2427 & 0.1462 \\
    \bottomrule
    \end{tabular}
  \label{tab:late_change1}
\end{table}

The tests were constructed as follows.  First, the epidemic is declared on
day $113$ with deployment under the Robust Policy beginning on day $140$,
i.e. $28$ days after the declaration.  As we saw in Example 1, the costliest tuple should the
Robust Policy be implemented is $(0.01168,0.0135,140)$ (thus, $p$ changes on the day of deployment).
Table \ref{tab:late_change1} evaluates the Robust Policy (and the No-Intervention policy) against this tuple, and also against tuples of the form
$(0.01168,\tilde{p},\tilde{d})$, with $\tilde{p} = 0.014, \, \mbox{or} \, 0.015$ and $\tilde{d} = 150, 155, 160,$ and $165$ and greater.
Table \ref{tab:late_change2} displays the same evaluations, but against the No-Action-Max-Cost tuple $(0.01092,0.0135,140)$, and against tuples of the form
$(0.01092,\tilde{p},\tilde{d})$, with $\tilde{p} = 0.014,\, \mbox{or} \,\, 0.015$ and $\tilde{d} = 170, 175, 180$ and greater.


\begin{table}[htbp]
  \centering
  \caption{Effects of a late change in $p$ on Scenario 1}
  \vspace{-1pt}
    \begin{tabular}{r|c|ccc|ccc}
    \multicolumn{ 8}{c}{{ Worst-case scenario given that no policy is implemented}} \\[3pt]
    \toprule
    {\bf Scenarios} & {\bf Original } & \multicolumn{ 3}{c}{{\bf Hypothetical 1}} & \multicolumn{ 3}{|c}{{\bf Hypothetical 2}} \\
    \midrule
    $p_1$   & 0.01092	& \multicolumn{ 3}{c}{0.01092} & \multicolumn{ 3}{|c}{0.01092} \\
    $p_2$  & 0.0135 & \multicolumn{ 3}{c}{0.014} & \multicolumn{ 3}{|c}{0.015} \\
    day epidemic declared & 133   & \multicolumn{ 3}{c}{133} & \multicolumn{ 3}{|c}{133} \\
    day deployment starts & 140   & \multicolumn{ 3}{c}{140} & \multicolumn{ 3}{|c}{140} \\
    day change $p$ & 140   & 170   & 175   & 180  & 170   & 175   & 180  \\
    cost Robust Policy & 0.0487 & 0     & 0     & 0     & 0.9827 & 0.2766 & 0 \\
    cost No Intervention & 4.58 & 0.7684 & 0.0093 & 0     & 2.7795 & 1.1467 & 0.0399 \\

    \bottomrule
    \end{tabular}
  \label{tab:late_change2}
\end{table}
In the case of Table \ref{tab:late_change1} we see that for $\tilde{d} = 150$ and  $\tilde{d} = 155$ the cost of the Robust Policy does increase over its worst-case cost under the
original robust model ($0.0508$) -- however the increase is modest relative
to the cost of not intervening.  Most importantly,  the cost of the
Robust Policy decreases rapidly
as $\tilde{d}$ moves out of scope of the original model.  In fact, the same is true even
under the No-Intervention Policy.

This is the salient fact
that we want to point out.  Its explanation is simple: even though $p$ is increasing to a larger value than originally envisioned, if $\tilde{d}$ is large, the impact of the increase is negligible -- the increase is taking place so late, that by then the epidemic (under the
initial $p$ value) has largely run its course.  In other words, we have a
long-running epidemic which, {\em under the Robust Policy}, ends up having no
impact.

Along the same lines, it is also worthwhile to note that in all cases where the out-of-sample increase in $p$
actually does have a cost impact, the change takes place shortly after the epidemic declaration and very shortly after deployment of the robust strategy begins. These facts indicate that there is a critical period (in this example of duration less than forty days or even shorter) starting from the declaration of the epidemic during which \textit{correct} action is important.  To some degree this supports our model of rolling-out a fixed strategy that deals with the immediate future; should the epidemic
``slow down'' only to restart much later, a completely new plan would then
be deployed.

\section{Extensions - policies with recourse}\label{thefuture}
Here we outline two extensions to our approach where the decision-maker can
apply recourse when conditions significantly deviate from predictions.  From
a broad optimization
perspective such strategies clearly make sense -- why stick with a rigid strategy?.  However, we stress that in view of the experiments at the end of the last section, the likely
benefit from a dynamic policy would be primarily realized during a period immediately
following the epidemic declaration and of relatively short duration. As as a result, a decision maker would likely
face severe logistical constraints if attempting to rapidly redeploy large
numbers of staff.  Thus, a dynamic policy might only be able to perform
\textit{small} changes on a day-to-day basis.  From an operational perspective such small
changes would of course make sense and would be applied; however the recourse might only
amount to a set of small adjustments to a pre-computed strategy.

In both cases that we consider we assume
that a decision-maker can only approximately observe the behavior
of $p_t$ (the value of $p$ at time $t$).  The Benders-decomposition solution methodology that we describe in the Appendix can be adapted to handle both models.

\subsubsection{Robust optimization with recourse}\label{robopt}
Assume that uncertainty of $p_t$ is modeled using $m$ intervals,
or $tranches$, denoted $I_1, I_2, \ldots I_m$, for some integer $m > 0$.
In addition, two time
periods are given, $J^{min}$ and $J^{max}$.  A realization of the
uncertain values $p_t$ plays out as follows:
\begin{itemize}
\item[(i)] A period $J^{min} \le J \le J^{max}$ is chosen, as well as a value $p^{(1)} \in I_1$,
a tranche $H$, and a value $p^{(2)} \in H$.
\item[(ii)] For each $t < J$ we have $p_t = p^{(1)}$ whereas for
each $t \ge J$, $p_t = p^{(2)}$.
\end{itemize}
\noindent In other words, $p_t$ is allowed to change values, once.  The
decision-maker operates as follows.  First, a fixed time period, $S$, in which
the deployment strategy will be \textit{revised} is
chosen in advance.  Then
\begin{itemize}
\item [(a)] At time $t = 1$,
the decision-maker knows that $p_1 \in I_1$, but
does not know its precise value. The decision-maker then produces an initial,
or ``first-stage''
deployment plan indicating the amount of staff
to call-up at time $t$, for each $t < S$.  Further,
$m$ alternative ``second-stage'' deployment plans (labeled $1, 2, \ldots, m$)
are announced,
each of which specifies the level of staff to deploy at each period $t \ge S$.
Each second-stage plan is compatible with the first-stage plan in terms of
deployment constraints (such as total staff availability).
\item [(b)]  At time $t = S$ the
decision-maker observes the tranche that $p_S$ belongs to, but not the actual
value of $p_S$.  If tranche $h$ is observed, then second-stage plan $h$
is rolled out on periods $t \ge S$.
\end{itemize}

\noindent \textit{Example.}  Consider a case with three tranches:
$I_2 = [.1,  .15], \ I_1 = [.15, .18], \ \mbox{and} \ I_3 = [.18, .19]$.
Further, $J^{min} = 18$ and  $J^{max} = 60$. A realization of
the data is the following: $J = 20$, $p^{(1)} = .155$, $H = I_3$ and $p^{(2)} = .19$.  Thus, for $t < 20$
we have $p_t = .155$ while for $t \ge 20$, $p_t = .19$.  Suppose $S = 30$.
The decision-maker knows that $.15 \le p_1 \le .18$; at time $S = 30$ the
decision-maker observes the tranche $I_3$ and thus knows that for $t \ge 30$,
$.18 \le p_t \le .19$, leading to an appropriate set of deployments
for future periods.\\

\noindent The model we just described is inherently adversarial: it assumes
that, subject to the stated rules, data will take on worst-case attributes.
This feature can, of course, lead to overly conservative plans.  However,
the proper statement to make is that the plans will be conservative only
if the data (the tranches, in particular) allow it to be so.  A different drawback in
this model is that an adversary could simply ``wait'' to change $p$ until just after
the review period
$S$.  As the experiments in Section \ref{outofsample} indicate, if $S$ is not too small the
impact of such a late change in $p$ is much decreased.  If $S$ is chosen too large, then
the adversary can change $p$ much earlier, and then the impact will be felt before the
review can take place.  Thus, care must be chosen in selecting $S$.  A review strategy
that proceeds dynamically is given next.

\subsubsection{Stochastic optimization}\label{stochopt}
In contrast to the above, we also outline a model where the contagion probabilities $p_t$ behave stochastically. Assume:
\begin{itemize}
\item [(i)] $p_1$ takes a random value drawn from a uniform distribution over
a $known$ interval $Z$.
\item [(ii)] For any $t$, $p_{t+1} = p_t \ + \ \epsilon_t$, where $\epsilon_t$
is normally distributed, with zero mean and small standard deviation (small
compared to the width of $Z$).
\end{itemize}
\noindent Thus, $p_t$ executes a random walk, starting from $Z$.
Let $\bar p$ indicate the center (mean value) of
$Z$.
The decision-maker's actions, in this model,
are as follows.

\begin{itemize}
\item[(a)] At time $t = 1$, the decision-maker produces a policy consisting of a
triple $(\tilde h, \, \alpha, \, R)$ where
$\alpha \ge 0$, and $\tilde h$ is an initial staff surge
deployment plan indicating for each $t$ the quantity $\tilde h_t$,
the amount of staff to call-up at time $t$.  Finally, the decision-maker also places
an additional quantity $R \ge 0$ of surge-staff on $reserve$ (i.e., not yet part of
the deployment plan).   For notational purposes, write $h_t = \tilde h_t$, for each $t$.
\item [(b)] The decision-maker can revise the plan at any of several
checkpoints $t_1 < t_2 < \ldots < t_r$.
\item [(c)] Consider checkpoint $t_i$ ($1 \le i \le m$). Let
\begin{eqnarray}
\Theta(t_i, \bar p) & = & \mbox{the total number individuals infected by time $t_i$, in the SEIR model where $p_t = \bar p$,} \nonumber \\
&& \mbox{for $1 \le t_i$} \nonumber \\
\Gamma(t_i) & = & \mbox{the $actual$ total number of individuals who are infected by time $t$.} \nonumber
\end{eqnarray}
Thus, $\Theta(t_i, \bar p)$
is an estimate of the random variable $\Gamma(t_i)$.  We assume that the decision-maker
observes, as an estimate for $\Gamma(t_i)$, a quantity
\begin{eqnarray}
\nu(i) & = & \Gamma(t_i) \ + \ \delta(t_i), \label{noisyobservation}
\end{eqnarray}
\noindent where $\delta(t_i)$ is normally distributed with zero mean and known
variance.  Then, the deployment plan is revised as follows. Write
$$ \Delta = \, \frac{\nu(i) \ - \ \Theta(t_i, \bar p)}{t_{i}}.$$

Then, the decision-maker resets
\begin{eqnarray}
h_t & \leftarrow & h_t \ + \ \alpha\,\Delta,  \ \ \ \mbox{for each $t_i \le t < t_{i+1}$, and} \label{update} \\
R & \leftarrow & R \ - \ (t_{i+1} - t_{i}) \, \alpha \, \Delta. \label{downdate}
\end{eqnarray}
\end{itemize}
\vskip -5pt
\noindent \textit{Comment.} The proposed scheme constitutes an example of
affine control.  The decision-maker computes an initial estimate of the
deployment plan (the initial $h_t$) which are corrected on the basis of
real-time observations, using the parameter $\alpha$ to modulate the
corrections.  The quantity $\Delta$ indicates the (per-period) deviation
between the observed behavior of the epidemic (including ``noise'', given
by $\delta$) and the predicted behavior as given by $\Theta(t_h, \bar p)$.

When $\Delta > 0$, the behavior of the epidemic is worse
than initially estimated and Rule
(\ref{update}) will increase the near-future surge staff deployment.
This increase is drawn from the surge staff held in reserve, as per
equation
(\ref{downdate}).  When $\Delta < 0$, (\ref{downdate}) moves staff back
into the reserve ranks.

Note that in order for Rule (\ref{downdate}) to be feasibly applied,
the resulting $R$ must be $nonnegative$.  For this to be the case
when $\Delta > 0$ we should
have $\alpha$ small enough, and, especially, that the initial estimate
for $R$ be ``high enough.''  These are constraints to be maintained in our
solution algorithms.

In order to consider algorithmic implications of this model, consider a fixed \textit{sample path} specifying
a value for each $p_t$ and for each $\delta(t_i)$ (see equation (\ref{noisyobservation})).  Denote by $T$ the length of the planning horizon. Recall that in
the stochastic model the decision-maker announces, at $t = 1$,
a triple $(\tilde h, \alpha, R)$.  We have:

\begin{LE} Under the given sample path
the cost incurred by the policy $(\tilde h, \alpha, R)$ (as per rules (b)-(c) )
is given by a linear program whose right-hand side is an affine function
of the $\tilde h_t$ and $\alpha$.   \end{LE}

\noindent As a consequence, we expect that several well-known methodologies for stochastic optimization, such as sample average approximation \cite{shapiro1}, \cite{shmoys1} and
stochastic gradient schemes \cite{robbinsmonro}, \cite{kushnerclark}, will be well-suited for our problem.

\section{Discussion}

The potential impact of a severe influenza epidemic on essential services
is a matter of critical importance from a Public Health perspective.
There are current concerns that an influenza virus might mutate into a highly contagious strain for which humans have little or no immunity, resulting in a rapid widespread of the disease. In the event of such an epidemic, or pandemic, organizations that provide critical social infrastructure, such as health care clinics, police departments, public utilities,
food markets and supply chains are likely to face severe workforce shortage that
would jeopardize their continuity of operations. In particular, at the most critical stage of the epidemic, health care clinics would be required to provide care for an extraordinary number of patients, and thus could ill-afford staff shortfalls. To address these issues, the CDC and the HHS Department have urged organizations to design contingency plans. As part these plans, organizations would deploy emergency 'surge' staff to compensate for a potential deficit in workers' availability.

The intelligent design of such contingency plans would necessarily model the
spread of the epidemic.  In the case of a new strain of the influenza virus, such modeling entails incorporating significant uncertainty: as far as we know, there are no clear ways of accurately predicting the transmissibility of a new virus strain. Moreover, it has been suggested that the effective contagion rate may change under different environmental conditions, such as weather, which are likewise unpredictable. Additionally, social patterns may change during the epidemic, for example due to the implementation of public health measures such as quarantine or social distancing.

In this paper we study how to design robust pre-planned surge staff deployment strategies that would help an organization cope with workforce shortfalls while optimally hedging against uncertainty. We propose fast and accurate algorithms that prove sufficiently flexible to incorporate intricate uncertainty sets that reflect the changeability of the environment. Our goal is to bring insights on the structure of optimal robust strategies and on practical rules-of-thumb that can be deployed during the epidemic.


There are many extensions to this work. In the case of health care workers, rates of transmission may be minimized during a severe pandemic by means of protective equipment or the use of antiviral medication.  Another extension could be to incorporate the use of cross-staffing, either within a hospital or a network of clinics and institutions. In this study we are considering the entire workforce of a particular organization as a whole; however, as part of the contingency plans proposed by the CDC, organizations are encouraged to identify critical staff positions and cross-train personnel accordingly. The model could also be furthered refine by modeling varying levels of exposure according to different subgroups of employees. Finally, from
an operational perspective, organizations may be able to react to evolutions of the
epidemic that significantly differ from initial estimations; how to incorporate
optimal recourse decisions is a topic we plan to address in the sequel.

\appendix
\section{Appendix - Solving the Robust Problem} \label{appendix}
Here we will describe an algorithm for solving the robust optimization problem (\ref{eq:V_star_def}).  The algorithm can be considered a variant of
Benders' decomposition \cite{benders}.  We begin by characterizing the representation of the cost incurred by a given deployment vector
under a given vector of contagion probabilities $\vec{p}$ as a linear program, a critical ingredient of our procedure.

\subsection{The robust problem as an infinite linear program}\label{subsec:SEIRsurgestaff}

Let $h$ be a given surge staff deployment vector, and let $\vec{p}$ describe
an evolution of the contagion probability.  We will next
obtain a representation
of the quantity $V(h \, | \, \vec{p})$ as the value of a linear program; this
will be a key step in the solution of problem (\ref{eq:V_star_def}).  To this
effect, let $h$ and $\vec{p}$ be given, let ($S^j_0, E^j_0, I^j_0, R^j_0, \; j = 1,2$) denote the initial $SEIR$ distribution; and
let all the other $SEIR$ parameters be fixed.
As noted before, this fully describes the predicted epidemic trajectory.
In particular, we obtain the values $\beta_t$ for all $t$,
which in turn fully define the disease dynamics of the procured staff
through equation (\ref{eq:SEIR_hires01}) given below.

The next step is to obtain a linear inequality representation of the quantity of available surge staff
at any given time $t$. Define:
\begin{itemize}
    \item $S^3_{t,j},$ the number of surge staff deployed at time $t - j$ that remain susceptible on day $t$, and
    \item $E^3_{t,j}$, the number off exposed surge staff deployed at time $t - j$.
\end{itemize}
As before, $T$ is the time horizon for the epidemic.  Using this notation, we have that:

\begin{eqnarray}\label{eq:SEIR_hires01}
    \begin{aligned}
        S^3_{t+1, 1} &= h_t, & t &= 1,..., T - \tau  \\
        S^3_{t+1, j+1} &= e^{-\lambda^2_t \beta_t p_t} \, S^3_{t,j}, & t &= 2,..., T - 1, \ \ \ \ j = 1,\ldots, \min\{t - 1, \tau - 1\}\\
        E^3_{t+1, 2} &= (1 - e^{-\lambda^2_t \beta_t p_t}) \, S^3_{t,1}, & t &= 2,...,T - 1 \\
        E^3_{t+1, j+1} &= (1 - e^{-\lambda^2_t \beta_t p_t}) \, S^3_{t,j} + e^{-\mu_{E_2}} E^3_{t,j}, & t &= 3,...,T - 1, \ \ \ \ j = 2,\ldots, \min\{t - 1, \tau - 1\}\\
    \end{aligned}
\end{eqnarray}
\noindent We note that the parameters $\lambda^2_t, \mu_{E_2}$ correspond to those of the workforce subgroup (and not the general population). This follows the assumption that the surge staff will be in the same circumstances as the workforce of interest.

Given \ref{eq:SEIR_hires01} and assuming that as soon as a staff member (regular or surge) becomes infectious he/she does not show up for work, the available \textit{surge} staff at each point in time $t$ is given by
\begin{equation}\label{eq:avextras_t}
    \sum_{j = 1}^{\min\{t,\tau\}}S^3_{t,j} + \sum_{j = 1}^{\min\{t,\tau\}}E^3_{t,j}, \ \ \mbox{for } t = 2, \ldots, T - \tau,
\end{equation}
where we have set $E^3_{2,j} = 0$.

Similarly, the total available \textit{regular} staff at time $t$ is given by the sum
\begin{equation}\label{eq:avreg_t}
    S^2_t + E^2_t + R^2_t.
\end{equation}

\noindent Thus, the sum of (\ref{eq:avreg_t}) and (\ref{eq:avextras_t}) yields the total available workforce at time $t$.

Now we return to the robust optimization problem (\ref{eq:V_star_def}).  For each period $t$, let $\omega_t$ denote the total number of available workers at time  $t$.
Based on discussions above, we assume that the cost function to be optimized is given by a convex piecewise-linear
function of the form $\max_{1 \le i \le L} \left\{ \sigma_i \, x_t + k_i \right\}$, for appropriate $L$, $\sigma$ and $k$.  This holds for the threshold
function case (section \ref{subsec:threshfun}) and for the queueing
case (section \ref{subsec:queue})
Using notation as in eq. (\ref{eq:costdef}), we therefore have
\begin{eqnarray}\label{eq:SEIR_Vhp}
    \begin{aligned}
    V(h \, | \, \vec{p}) &:=  \min \,\, \sum_{t = 1}^T z_t \\
    \mbox{s.t.} && \\
        &S^3_{t+1, 1} = h_t,  &t& = 1,..., T - \tau  \\
        &S^3_{t+1, j+1} = e^{-\lambda^2_t \beta_t p_t} \, S^3_{t,j},  &t& = 2,..., T - 1, \, j = 1,\ldots, \min\{t - 1, \tau - 1\}\\
        &E^3_{t+1, 2} = (1 - e^{-\lambda^2_t \beta_t p_t}) \, S^3_{t,1}, &t& = 2,...,T - 1 \\
        &E^3_{t+1, j+1} = (1 - e^{-\lambda^2_t \beta_t p_t}) \, S^3_{t,j} + e^{-\mu_{E_2}} E^3_{t,j}, &t& = 3,...,T - 1, \, j = 2,\ldots, \min\{t - 1, \tau - 1\}\\
        &\omega_t = S^2_t + E^2_t + R^2_t + \sum_{j = 1}^{\min\{t-1,\tau\}}S^3_{t,j} + \sum_{j = 1}^{\min\{t-1,\tau\}}E^3_{t,j}, &t& = 1,...,T\\
        & z_t \, \geq \sigma_i \omega_t + k_i,   &1 &\le i \le L, \, \,\, t \, = 1,...,T
    \end{aligned}
\end{eqnarray}

\noindent In this formulation $S$, $E$, $\omega$ and $z$ are the variables (they should be indexed by $\vec{p}$, but we omit this for simplicity of notation).
Notice that the quantity $h_t$ appears explicitly in only the first set of constraints. Also, by construction all variables are
nonnegative.  We can summarize this formulation in a more compact form.  Using appropriate matrices $A_{\vec{p}}$ and $C_{\vec{p}}$ and vectors $\kappa_{\vec{p}}$ and $d_{\vec{p}}$,
\begin{eqnarray}
       V(h | \vec{p}) \ := \ && \, \underset{x}\min \, \kappa_{\vec{p}}^T x  \label{eq:LP_p_1} \\
                    \mbox{s.t.} && A_{\vec{p}}\, x \, = \, h    \label{eq:LP_p_2}\\
                    && C_{\vec{p}} \, x \, \geq d_{\vec{p}}            \label{eq:LP_p_3}
\end{eqnarray}
\noindent Here, $x$ denotes a vector of auxiliary variables comprising all
$S$, $E$, $w$ and $z$.    We can now write:

\begin{eqnarray}\label{eq:LP_p}
    \begin{aligned}
        V^* \ :=
               \ \  &  \underset{h, x}\inf & \nu \\
                & \mbox{s.t.} & \nu &\geq \kappa_{\vec{p}}^T x,  & \forall \vec{p} \in \cP  \\
                & & A_{\vec{p}} \, x &= h,  & \forall \vec{p} \in \cP  \\
                & & C_{\vec{p}} \, x &\geq d_{\vec{p}},         & \forall \vec{p} \in \cP  \\
                & & h &\in \cH.
    \end{aligned}
\end{eqnarray}
\noindent Problem (\ref{eq:LP_p}) is an \textit{infinite linear program}. It
can be shown that the ``inf'' is a ``min''.

\subsection{Algorithm} \label{bendersdecomposition}
We can now describe our procedure for solving
optimization problem problem (\ref{eq:LP_p}).

\vspace{.1in}

\noindent Let $\pi$, $\alpha$ be \textit{feasible} dual vectors for the LP (\ref{eq:LP_p_1})-(\ref{eq:LP_p_3}), corresponding to (\ref{eq:LP_p_2}) and (\ref{eq:LP_p_3}), respectively. Then by
weak duality,
\begin{equation}\label{eq:weak_dual}
    V (h \, | \, \vec{p}) \geq \pi^T h + \alpha^T d_{\vec{p}}.
\end{equation}
Furthermore, $\pi$ and $\alpha$ are optimal for the dual if and only if
\begin{equation}
    V (h \, | \, \vec{p}) = \pi^T h + \alpha^T d_{\vec{p}}.
\end{equation}
Denoting by $D(\vec{p})$ the set of feasible duals to LP (\ref{eq:LP_p_1})-(\ref{eq:LP_p_3}), it follows that
\begin{equation}
    V (h \, | \, \vec{p}) = \underset{(\pi,\alpha) \in D(\vec{p})} \max \pi^T h + \alpha^T d_{\vec{p}}.
\end{equation}
and so we can rewrite (\ref{eq:V_star_def}) as
\begin{equation}\label{eq:V_star_minmaxmax}
    V^* =  \underset{h \in H}{\min} \; \underset{\vec{p} \in P}{\max} \; \underset{(\alpha, \pi) \in D(\vec{p})}{\max} \alpha^Th + \pi^Td_{\vec{p}}.
\end{equation}

\noindent Now consider a finite family of vectors
$(\pi_k, \alpha_k)$ ($k = 1, \ldots, K$) such that for each $k$ there is a vector $\vec{p}(k) \in \cP$ with
$(\pi_k, \alpha_k) \in D(p^k)$.
Then, by (\ref{eq:V_star_minmaxmax}), we have
\begin{equation}\label{eq:V_lower}
    V^* \geq  \underset{h \in H}{\min} \; \underset{1 \le k \le K}{\max} \; \alpha^T_k h + \pi^T_k d_{\vec{p}(k)}.
\end{equation}
\noindent This observation gives rise to the following algorithm.

\begin{center}
  \textsc{{\bf Algorithm B}}\vspace*{5pt}\\
  \fbox{
    \begin{minipage}{0.95\linewidth}
\hspace*{.1in} \\
      \hspace*{.1in} {\bf 0.} Set $K = 0$ and $r = 1$.\\ \\
      \hspace*{.1in} {\bf 1.} Let $h^r$ be an optimal solution for the LP
  \begin{eqnarray*}
    \begin{aligned}
        W^r :=  \, &\underset{z,h}{\min} && z \\
                &\mbox{s.t.} & &  z \geq \alpha_{k}^Th + \pi_{k}^Td_{\vec{p}(k)}, \;\; 1 \le k \le r-1 \label{eq:cutplane}\\
                & & & h \in \cH.
    \end{aligned}
  \end{eqnarray*}
    \hspace*{.1in} {\bf 2.} Let $\vec{p}(r) \in \cP$ be such that
        \begin{equation}\label{eq:oracle}
        V(h^r \, | \, \vec{p}(r)) = \underset{ \vec{p} \in \cP} \max \, V (h^r\, | \, \vec{p}),
        \end{equation}

        \hspace*{.3in} \textit{Note}: $\vec{p}(r)$ is the contagion probability vector that attains the worst case should deployment\\  \hspace*{.3in} vector $h$ be used. \\

        \hspace*{.1in} {\bf 3.} Using notation as in formulation (\ref{eq:LP_p_1})-(\ref{eq:LP_p_3}), let $(\pi_{r}, \alpha_{r})$ be optimal duals for the LP
\begin{eqnarray}
       \ && \, \underset{x}\min \, \kappa_{\vec{p}(r)}^T x  \label{repeateq:LP_p_1} \\
                    \mbox{s.t.} && A_{\vec{p}(r)}\, x \, = \, h^r    \label{repeateq:LP_p_2}\\
                    && C_{\vec{p}(r)} \, x \, \geq d_{\vec{p}(r)}            \label{repeateq:LP_p_3}
\end{eqnarray}
        \hspace*{.3in} \textit{Note}: The value of this LP is $V(h^r | \vec{p}(r))$.\\

    \hspace*{.1in} {\bf 4.} {\bf If} $W^r \ge V(h^r \, | \, \vec{p}(r))$, \hspace*{.1in} {\bf STOP} -- algorithm has terminated. \\

 \hspace*{.3in} {\bf Else}, reset $r \leftarrow  r + 1$ and go to {\bf 1}.\\

\end{minipage}
  }
\end{center}

\noindent \textit{Remark}: the above algorithm can be viewed as a special case of Benders' decomposition \cite{benders}.  The linear program solved in Step 1 is called
the \textit{master problem}.

\begin{LE}\label{le:benders} For any $r \ge 0$ we have (a) $W^r \leq W^{r+1}$ and (b)
$W^r \leq V^* \leq V(h^r\, | \, \vec{p}(r))$.
\end{LE}
\noindent \textit{ Proof.} (a) This follows because in each iteration we add a
new constraint to the problem solved in Step 1. (b) This follows as
per equation (\ref{eq:V_lower}). \QED

\begin{CO} If the algorithm terminates in Step 4, it has correctly
solved problem (\ref{eq:V_star_def}). \end{CO}

\noindent In practice we would not stop the algorithm in Step 4, but rather use part (b)
of Lemma \ref{le:benders} to stop when a desired optimality guarantee is reached.  An issue of theoretical interest is the rate of convergence of Algorithm B. The following result indirectly addresses this question.

\begin{LE}\label{le:ellipsoid} There is an algorithm that computes $V^*$ by
solving problem (\ref{eq:oracle}) a polynomial number of times.\end{LE}

\noindent We stress that the algorithm in Lemma \ref{le:ellipsoid} is \textit{not}
Algorithm B; rather, it relies on the well-known equivalence between separation
and optimization \cite{ellipsoid} and it is similar to Algorithm B except that
(essentially) the computation of $h^r$ in Step 1 is performed differently.
The resulting algorithm, while theoretically efficient, may not be practical.
Instead, in practice researchers in the nonconvex optimization community
would rely on a classical cutting-plane algorithm such as Algorithm B.

This discussion highlights the importance, both theoretical and practical,
 of Step 2 of the algorithm, namely the computation of
$\argmax_{ \vec{p} \in \cP} ~ V(h \, | \, \vec{p})$ for a given deployment vector $h$.  This is a one-dimensional maximization problem,
but possibly non-concave (recall Figure \ref{chart:C_p_3H}).  In this context, an important experimental observation is that
the computation of $V(h \, | \, \vec{p})$ for given $h$ and $\vec{p}$ is extremely fast, even for large $T$ -- typically,
$V(h \, | \, \vec{p})$ can be computed in approximately one {\it hundred thousandth} of a second
on a modern computer.  In our implementation we make use of this observation by discretizing the
set $\cP$.

Recall that in our uncertainty model (Section \ref{subsec:UncertModels}), that
for $t \le \breve t$, $p_t$ takes a fixed value in $[p^L, p^U]$ whereas
for $t > \breve t$, $p_t$ takes a fixed value in $[\hat p^L, \hat p^U]$.  Let $N$ be a large integer, and write
\begin{eqnarray}
&& \cQ_N \ = \ \left\{ \, p^L + \frac{p^U - p^L}{N} \, j, \ 0 \le j \le N \right\} \ \ \mbox{and} \ \ \hat \cQ_N \ = \ \left\{ \, \hat p^L + \frac{\hat p^U - \hat p^L}{N} \, j, \ 0 \le j \le N \right\}
\end{eqnarray}
\noindent $\cQ_N$ is a finite approximation to the interval $[p^L, p^U]$  and likewise with  $\hat \cQ_N$ and $[\hat p^L, \hat p^U]$.  Define
$$\Pi_N \ \ = \ \ \{ \, (q, \hat q, \breve t) \ : \ q \in \cQ_N, \ 1 \le \breve t \le T, \ \hat q \in \hat \cQ_N \}, $$
\noindent and for any $\pi = (q, \hat q, \breve t) \in \Pi_N$,  let $\vec{p}(\pi)$ be the vector of probabilities defined by
$$p_t(\pi) \, = \, q \ \ \mbox{for $1 \le t \le \breve t$}, \ \ \mbox{and} \ \ p_t(\pi) \, = \, \hat q \ \ \mbox{for $\breve t < t \le T$}.$$
Finally, define
$$\cP_N \ \ = \ \ \{ \ \vec{p}(\pi) \ \ : \ \ \pi \ \in \ \Pi_N \ \}.$$
Thus, $\cP_N$ is a discrete approximation to $\cP$, and, furthermore
$V^*_N  :=  \min_{ h \in \cH} \max_{ \vec{p} \in \cP_N} ~ V(h \, | \, \vec{p})$
can be computed using a (finite) linear program analogous to (\ref{eq:LP_p})
which we include for completeness:
\begin{eqnarray}\label{eq:LP_pN}
    \begin{aligned}
        V^*_N :=
                & \underset{h, x}\min & \nu \\
                & \mbox{s.t.} & \nu &\geq \kappa_{\vec{p}}^T x^{\vec{p}},  & \forall \vec{p} \in \cP_N  \\
                & & A_{\vec{p}} x^{\vec{p}} &= h,  & \forall \vec{p} \in \cP_N  \\
                & & C_{\vec{p}} x^{\vec{p}} &\geq d_{\vec{p}},         & \forall \vec{p} \in \cP_N  \\
                & & h &\in \cH.
    \end{aligned}
\end{eqnarray}
Clearly, $V^*_N \le V^*$ since $\cP_N \subset \cP$, and one can show that $V^*_N \rightarrow V^*$
as $N \rightarrow +\infty$.  From a practical perspective, a large enough
value for $N$, such as $N = 1000$, should provide an excellent approximation
to our robust staffing problem. Note that, for any $N$ we have $\cP_N \subset \cP_{2 N}$.  Hence, we can proceed by computing $V^*_1$, followed by $V^*_2$,
followed by $V^*_4$, and so on until the desired accuracy (in terms of $p$) is
attained in a logarithmic number of iterations.  Notice that in this framework
the cutting planes (\ref{eq:cutplane}) discovered in each iteration of Step 1
performed when computing each value $V^*_N$ remain valid for the computation
of $V^*_{2 N}$ (precisely because $\cP_N \subset \cP_{2 N}$) and thus each
iteration is warm-started by the preceding iteration.  Additional speed-up
strategies are possible.

\subsection{Improved algorithm}\label{subsec:algo2}

While the above algorithm described is efficient, it may sometimes require many iterations to achieve a desirable tolerance.  We improve on our algorithm by
first (carefully) enumerating a small subset of tuples $\hat \Pi$ from our
uncertainty set, and introducing into the master the constraints
used to describe SEIR model and the corresponding objective function inequalities \ref{eq:LP_pN} for each $\vec{p} \in \hat \Pi$.

Inspired by the algorithm itself, we choose these tuples iteratively by finding the worst-case behavior for a given feasible deployment strategy, adding the corresponding constraints to the master, and resolving. Formally, this
approach is described as follows.

\begin{center}
  \textsc{{\bf Procedure A}}\vspace*{5pt}\\
  \fbox{
    \begin{minipage}{0.95\linewidth}
\hspace*{.1in} \\
      \hspace*{.1in} {\bf 0.} Set $Q_0 = \emptyset$ and $r = 1$.\\

      \hspace*{.1in} {\bf 1.} Let $\vec{p}(r) \in \cP$ be the tuple such that
        \begin{equation}
        V(h^r \, | \, \vec{p}(r)) = \underset{ \vec{p} \in \cP} \max \, V (h^r\, | \, \vec{p}).
        \end{equation}

      \hspace*{.1in} Let $\tilde{V}^r = \min\{ V(h^k \, | \, \vec{p}(k)) \}_{k = 1}^{r}$ and update $Q_0 \leftarrow Q_0 \cup \{\vec{p}(r)\}.$
        \\

      \hspace*{.1in} {\bf 2.} Let $h^r$ be an optimal solution for the LP
    \begin{eqnarray*}
    \begin{aligned}
        \tilde{W}^r :=
                & \underset{h, x}\min & \nu \\
                & \mbox{s.t.} & \nu &\geq \kappa_{\vec{p}}^T x^{\vec{p}},  & \forall \vec{p} \in Q_0  \\
                & & A_{\vec{p}} x^{\vec{p}} &= h,  & \forall \vec{p} \in Q_0  \\
                & & C_{\vec{p}} x^{\vec{p}} &\geq d_{\vec{p}},         & \forall \vec{p} \in Q_0  \\
                & & h &\in \cH \\
                & & x \geq 0.
    \end{aligned}
\end{eqnarray*}

    \hspace*{.1in} {\bf 3.} {\bf If} $|\tilde{V}^r - \tilde{W}^r| < \epsilon \tilde{W}^r$ or $r > K$, \hspace*{.1in} {\bf STOP} -- algorithm has terminated. \\
    \hspace*{.3in} {\bf Else}, reset $r \leftarrow  r + 1$ and go to {\bf 1}.\\

\end{minipage}
  }
\end{center}

Parameters $\epsilon$ and $K$ are fixed a priori; they represent a duality gap tolerance and a maximum number of iterations for Procedure A, respectively. It is worth mentioning that with each iteration the program acquires a significant number of constraints: for a single tuple, under the queueing scenario, there are close to 2,800 constraints and more than 2,500 variables (after pre-solving). Thus, we expect $K$ to be small (not larger than $10$).

Once Procedure A terminates, we start Algorithm B, with the proviso that
the Master Problem is initialized as:

\begin{eqnarray*}
    \begin{aligned}
        W^r :=  \, &\underset{z,h,x}{\min} && z \\
                &\mbox{s.t.} &  z & \geq \alpha_{k}^Th + \pi_{k}^Td_{\vec{p}(k)}, && 1 \leq k \leq r-1 \label{eq:cutplane2}\\
                &  & z &\geq \kappa_{\vec{p}}^T x^{\vec{p}},  && \forall \vec{p} \in Q_0  \\
                & & A_{\vec{p}} x^{\vec{p}} &= h, &&  \forall \vec{p} \in Q_0  \\
                & & C_{\vec{p}} x^{\vec{p}} &\geq d_{\vec{p}},  && \forall \vec{p} \in Q_0  \\
                & & & h \in \cH \\
                & & & x \geq 0.
    \end{aligned}
\end{eqnarray*}

Algorithm B is then run as before. At each iteration $r$, Step 2 discovers a
scenario $\vec{p}(r)$, and  Step 3 produces a dual vector $(\pi_r, \alpha_r)$
which gives rise to the cut $z \geq \alpha_{r}^Th + \pi_{r}^Td_{\vec{p}(r)}$ that is added to the master.  The effect of the
enhancement provided by Procedure A is, typically, to drastically shortcut the
number of iterations required by Algorithm B; essentially, the algorithm has
been hot-started with a very good initial representation of the critical
constraints needed to define problem (\ref{eq:V_star_def}).

To illustrate the performance of the enhanced Algorithm B, Table \ref{tab:RunningTimes_Ex1} presents a number of statistics in the case of Example 1 (Section \ref{subsubsec:Health_Ex1}).  The table shows the worst tuple $(p_1, p_2, d)$,
the lower and upper bounds on the robust optimization problem,
and the CPU time per iteration (in seconds) used to compute the worst tuple $(p_1, p_2, d)$ and to solve the corresponding linear programs with AMPL\footnote{Time is measured using \texttt{\_total\_solve\_elapsed\_time}, which reflects the elapsed seconds used by the solve commands.}.
The computation of the worst tuple is the most time-consuming task at each iteration.  At present we have implemented this step as a search process which could
be improved in a number of ways.  However, the algorithm appears effective;
note that in this case the total CPU time does not exceed $35$ seconds.

Table \ref{tab:RunningTimes_Ex1} also shows how fast the algorithm converged for this example. In terms of the enhancement, we chose the maximum number of iterations to be $K = 10$ and $\epsilon = 0.05.$ The enhancement was used for the first 8 iterations, after which only one more cut was added to the Master Problem. The resulting duality gap is slightly larger than $0.005\%.$

\begin{table}[htbp]
  \centering
  \small
  \caption{Bounds and CPU times per iteration for Example 1 (Section \ref{subsubsec:Health_Ex1}). The algorithm used the enhancement for the first 8 iterations and switched to our algorithm for one more iteration.}
    \begin{tabular}{rccccccc}
    \addlinespace
    \toprule
          & \multicolumn{ 3}{c}{{\bf Worst tuple}} & \multicolumn{ 2}{c}{{\bf Convergence Bounds}} & \multicolumn{ 2}{c}{{\bf CPU time (seconds)}} \\
    \midrule
    {\bf Iter} & {\bf p\_1} & {\bf p\_2} & {\bf day ch} & {\bf Lower} & {\bf Upper} & {\bf worst tuple} & {\bf AMPL} \\
    Enh\_1 & 0.01092 & 0.0135 & 140   & 0     & 4.581151 & 3.406 & 0.328 \\
    Enh\_2 & 0.01172 & 0.0135 & 140   & 0     & 0.710181 & 3.453 & 0.344 \\
    Enh\_3 & 0.01132 & 0.0135 & 140   & 0.007060 & 0.217431 & 3.453 & 0.343 \\
    Enh\_4 & 0.01081 & 0.0135 & 142   & 0.031251 & 0.132066 & 3.468 & 0.391 \\
    Enh\_5 & 0.01185 & 0.0135 & 140   & 0.048387 & 0.132367 & 3.453 & 0.391 \\
    Enh\_6 & 0.01168 & 0.0135 & 140   & 0.050479 & 0.073479 & 3.469 & 0.453 \\
    Enh\_7 & 0.01177 & 0.0135 & 140   & 0.050686 & 0.056242 & 3.562 & 0.422 \\
    Enh\_8 & 0.01111 & 0.0135 & 140   & 0.050686 & 0.052368 & 3.453 & 0.453 \\
    Bend\_1 & 0.01172 & 0.0135 & 140   & 0.050765 & 0.050768 & 3.437 & 0.453 \\
    \bottomrule
    \end{tabular}
  \label{tab:RunningTimes_Ex1}
\end{table}

\bibliographystyle{plain}
\bibliography{bibfile}


\tiny Thu.Jul.30.090524.2015@littleboy

\end{document}